%% file: main.tex
\documentclass[hidelinks,onefignum,onetabnum]{siamart251216}

\usepackage[utf8]{inputenc}
\usepackage{amsmath, amssymb, mathtools} 
\usepackage{hyperref} 
\usepackage{subcaption}
\usepackage{multicol}
\usepackage{multirow}
\usepackage{mymacros}
\usepackage{xcolor}
\usepackage[labelfont=bf]{caption}
\usepackage{pdflscape}
\usepackage[noend]{algpseudocode}
\usepackage{cleveref}

\usepackage{algorithm}
\usepackage{booktabs}
\usepackage{adjustbox}
\usepackage{etoolbox}
\usepackage[italicdiff]{physics}
\usepackage{todonotes}
\usepackage{rotating, graphicx}
\setuptodonotes{inline}
\newcommand{\gp}{{\mathcal{GP}}}
\newcommand{\bKKrr}{\mathbf{K}_{11}}
\newcommand{\bKKnr}{\mathbf{K}_{21}}
\newcommand{\bKrn}{\mathbf{K}_{12}}

\newcommand{\bKKrrinv}{\mathbf{K}_{11}^{-1}}
\newcommand{\ld}{\ensuremath{\text{logdet}}}
\newcommand{\inv}{{-1}}
\newcommand{\dK}{\ensuremath{\mathbf{dK}}} %identity matrix
\newcommand{\sI}{\ensuremath{\sigma^2\mathbf{I}}} %identity matrix with noise
\newcommand{\tKK}{\widetilde{\bK}}
\newcommand{\had}{{\textstyle\odot}}

\newcommand{\AnchorNet}{\ensuremath{\text{\textproc{AnchorNet}}}}

\newcommand{\myalgsheader}[0]

\Crefname{appendix}{Appendix}{Appendices}
\crefname{appendix}{appendix}{appendices}

\AddToHook{cmd/appendix/before}{%
  \crefalias{section}{appendix}%
}

\algnewcommand{\IIf}[1]{\State\algorithmicif\ #1\ \algorithmicthen}
\algnewcommand{\EndIIf}{\unskip\ \algorithmicend\ \algorithmicif}
\algnewcommand{\IElse}[1]{\State\algorithmicelse\ #1\ \algorithmicthen}
\algnewcommand{\IfThenElse}[3]{% \IfThenElse{<if>}{<then>}{<else>}
	\State \algorithmicif\ #1\ \algorithmicthen\ #2\ \algorithmicelse\ #3}

\ifpdf
  \DeclareGraphicsExtensions{.eps,.pdf,.png,.jpg}
\else
  \DeclareGraphicsExtensions{.eps}
\fi

\newsiamremark{remark}{Remark}
\newsiamremark{hypothesis}{Hypothesis}
\crefname{hypothesis}{Hypothesis}{Hypotheses}
\newsiamthm{claim}{Claim}
\newsiamremark{fact}{Fact}
\crefname{fact}{Fact}{Facts}

\headers{LP Training for GPs}{Poguntke et. al}

\title{Low-rank low precision Training of Gaussian Processes}

\author{Moritz Poguntke\thanks{Department of Mathematics, TU Chemnitz, Chemnitz, Germany(\email{moritz.poguntke@mathematik.tu-chemnitz.de},\email{martin.stoll@mathematik.tu-chemnitz.de}).}
\and Eda Oktay\thanks{Max Planck Institute for the Dynamics of Complex Technical Systems, Magdeburg, Germany (\email{oktay@mpi-magdeburg.mpg.de}).}
\and Yuanzhe Xi\thanks{Department of Mathematics, Emory University, USA (\email{yuanzhe.xi@emory.edu}). The work of Y. Xi is supported by the National Science Foundation
under award DMS-2338904.}
\and Erin Carson\thanks{Department of Numerical Mathematics,
Faculty of Mathematics and Physics,
Charles University(\email{carson@karlin.mff.cuni.cz}). Supported by the European Union (ERC, inEXASCALE, 101075632) and the Charles University Research Centre program No. UNCE/24/SCI/005.}
\and Martin Stoll\footnotemark[1]}

\begin{document}
\maketitle

% \title{Mixed Precision Training of Gaussian Processes}
	
% \author{Eda Oktay\footnote{Department of Numerical Mathematics, Faculty of Mathematics and Physics, Charles University, \{oktay, carson\}@karlin.mff.cuni.cz. We acknowledge funding from Charles University Research Program No. UNCE/SCI/023, GAUK project No. 202722, the Exascale Computing Project (17-SC-20-SC), a collaborative effort of the U.S. Department of Energy Office of Science and the National Nuclear Security Administration, and by the European Union (ERC, inEXASCALE, 101075632). Views and opinions expressed are those of the
% authors only and do not necessarily reflect those of the European
% Union or the European Research Council. Neither the European Union nor
% the granting authority can be held responsible for them.},\hspace{1pt} Moritz Poguntke \footnotemark[1], \hspace{1pt} Yuanzhe Xi\footnotemark[1]   \hspace{1pt} Erin Carson\footnotemark[1] \hspace{1pt} and Martin Stoll\footnote{Department of Mathematics, TU Chemnitz, martin.stoll@mathematik.tu-chemnitz.de.}}
% \date{}
% \maketitle

% REQUIRED

\begin{abstract}
Training Gaussian processes (GPs) at large scale is computationally
challenging because evaluating the log-marginal likelihood and its
gradients involves large, dense kernel matrices. We develop a scalable
low-rank GP training approach that combines \AnchorNet\ landmark
selection with a fully independent training conditional (FITC)
correction. The resulting objective and its analytic gradients can be
evaluated using low-rank factors, diagonal quantities, and small dense
matrices, without forming the full kernel matrix or its derivatives.
The complete training procedure is implemented in single precision and
evaluated on nine large-scale regression datasets. Across these
problems, the method reaches the reference predictive accuracy while
substantially reducing training time, although the required
approximation rank varies across datasets. We also study the numerical
behavior of the full-kernel predictive solve and show that single
precision can become limiting for ill-conditioned problems. Finally,
we investigate further reducing selected matrix operations to fp16 or
TF32. These lower-precision variants provide only modest additional
speedups and can become unreliable when the reduced system is poorly
conditioned.
\end{abstract}

% REQUIRED

% REQUIRED

\section{Introduction}
\label{sec:intro}

Surrogate models can serve as accurate and computationally efficient
replacements for expensive models or can represent an underlying
unknown function from observed data \cite{gramacy2020surrogates}.
Gaussian processes (GPs) \cite{williams2006gaussian} are widely used
for this purpose in Bayesian optimization and machine learning because
they provide both predictions at unseen inputs and a principled
quantification of predictive uncertainty.

The performance of a GP depends strongly on the choice of covariance
or kernel function and its associated hyperparameters
\cite{scholkopf2001learning}. A central task in GP training is therefore
to estimate these hyperparameters from the available training data
\cite{RaWi}. Both hyperparameter training and posterior prediction
require matrix computations involving a large, dense covariance
matrix. A direct implementation based on matrix factorization has
cubic complexity in the number of training points, which quickly
becomes prohibitive for large datasets. Consequently, considerable
effort has been devoted to structured and low-rank approximations that
reduce these computational costs
\cite{williams2006gaussian,alder2024energy,wu2024large,
dong2017scalable,quinonero2005unifying,Wagner2027}.

Further computational gains are possible by exploiting reduced-precision
arithmetic on modern hardware, where single- and half-precision
operations can provide substantially higher throughput and lower
memory traffic than double precision. This has motivated increasing
interest in low-precision algorithms for large-scale GP computations
\cite{maddox2022low}. At the same time, reduced precision introduces
additional numerical challenges because GP kernel matrices can become
severely ill-conditioned during hyperparameter optimization.

The contributions of this paper are two-fold. First, we develop a
scalable low-rank GP training approach based on landmarks selected by
the \AnchorNet\ algorithm \cite{cai2022fast} together with a fully
independent training conditional (FITC) correction. We derive the
corresponding objective and analytic gradients so that both the
objective and its derivatives can be evaluated using low-rank factors
and small dense matrices, without forming the full kernel matrix or
its derivatives. Second, we study the numerical behavior of this
training procedure in reduced precision on modern GPU architectures.
We demonstrate that the complete low-rank training procedure can be
carried out effectively in single precision and investigate the extent
to which selected matrix operations can be further reduced to half
precision or TF32. Our experiments also identify the conditioning
regimes in which these lower-precision computations become unreliable.

The remainder of the paper is organized as follows.
In \Cref{sec:gp}, we review the GP formulation and the kernel models
used in this work. In \Cref{sec:training}, we derive the low-rank FITC
objective and its gradients. \Cref{sec::anchornet} describes the
\AnchorNet\ landmark-selection procedure, and \Cref{sec:alg}
summarizes the resulting training algorithm.
\Cref{sec:experiments_new} presents the numerical experiments,
including the effects of the approximation rank and the numerical
behavior of the predictive solve. Finally, \Cref{sec:outlook_new}
investigates the use and limitations of precision below fp32.

\section{Gaussian Processes}
\label{sec:gp}

A Gaussian process (GP) is a collection of random variables such that
every finite subset has a joint Gaussian distribution
\cite{williams2006gaussian}. It is completely characterized by its
mean function
\[
m(x) := \mathbb{E}[f(x)]
\]
and covariance, or kernel, function
\[
k(x,x')
:=
\mathbb{E}\bigl[(f(x)-m(x))(f(x')-m(x'))\bigr].
\]
We write
\[
f\sim\gp(m,k)
\]
and assume a zero mean, $m(x)=0$, throughout this work.

Let $\mathbf X_{\mathrm{train}}\in\R^{n\times d}$ contain the training
inputs $x_1,\dots,x_n$ as its rows, and define the corresponding
latent function values by
\[
\mathbf f_{\mathrm{train}}
=
[f(x_1),\dots,f(x_n)]^T.
\]
We observe
\[
\mathbf y_{\mathrm{train}}
=
\mathbf f_{\mathrm{train}}+\boldsymbol{\varepsilon},
\qquad
\boldsymbol{\varepsilon}
\sim
\mathcal N(\mathbf 0,\sigma^2\mathbf I_n),
\]
where $\sigma>0$ and the observation noise is independent of the
latent process.

\begin{figure}
    \centering
    \begin{subfigure}{0.5\textwidth}
        \includegraphics[width=\textwidth]{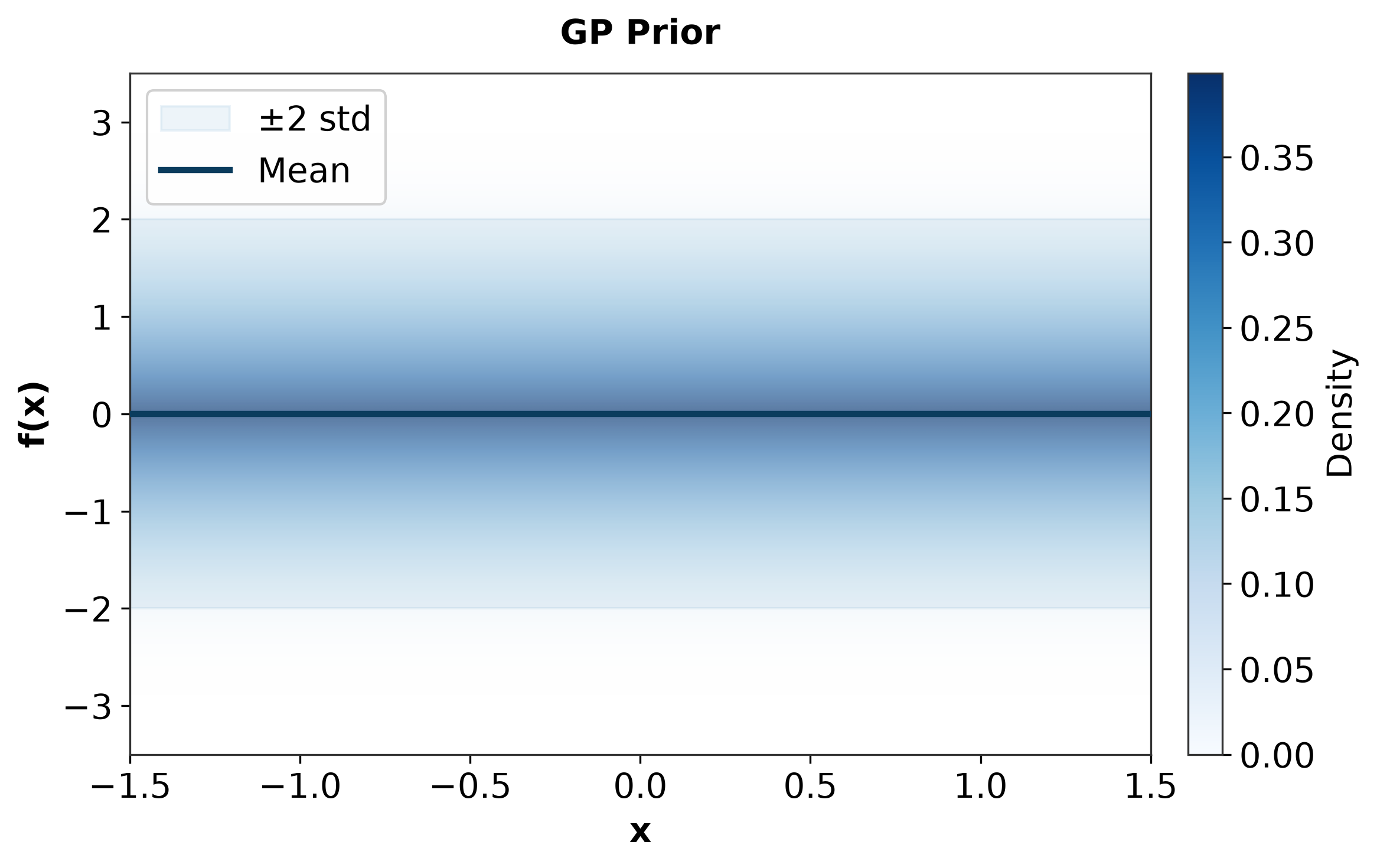}
        \caption{Prior with zero mean}
        \label{fig:prior}
    \end{subfigure}%
    \begin{subfigure}{0.5\textwidth}
        \includegraphics[width=\textwidth]{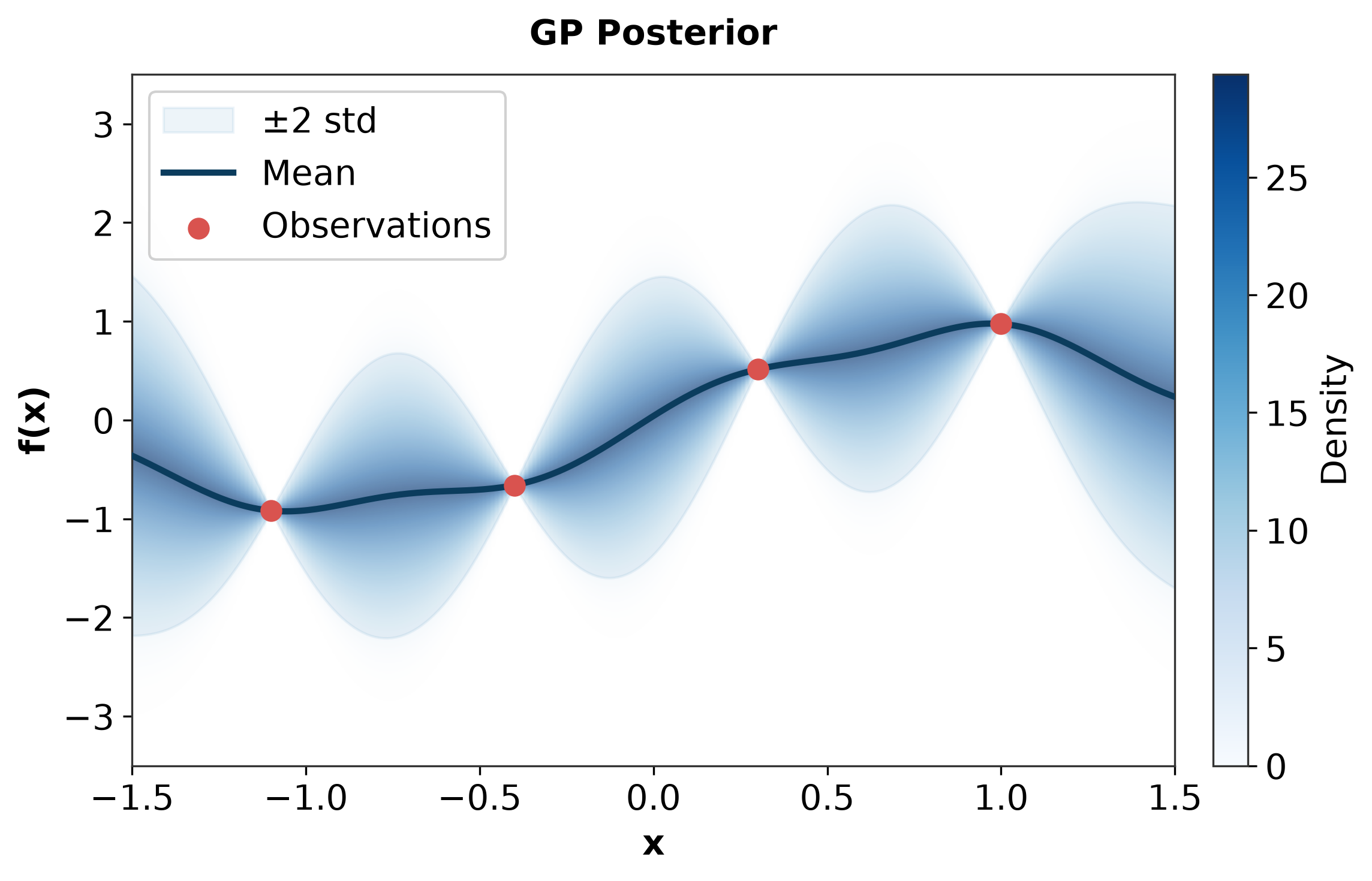}
        \caption{Posterior based on four data points}
        \label{fig:posterior}
    \end{subfigure}%
    \caption{A one-dimensional GP before and after conditioning on
    training data.}
    \label{fig:gp}
\end{figure}

For $n_*$ test inputs
$\mathbf X_{\mathrm{test}}\in\R^{n_*\times d}$, define
\[
\bK
=
k(\mathbf X_{\mathrm{train}},\mathbf X_{\mathrm{train}})
\in\R^{n\times n},
\]
\[
\bK_*
=
k(\mathbf X_{\mathrm{train}},\mathbf X_{\mathrm{test}})
\in\R^{n\times n_*},
\qquad
\bK_{**}
=
k(\mathbf X_{\mathrm{test}},\mathbf X_{\mathrm{test}})
\in\R^{n_*\times n_*},
\]
where each kernel matrix is formed by evaluating $k$ at the
corresponding pairs of inputs. Let
$\mathbf f_{\mathrm{test}}\in\R^{n_*}$ collect the latent function
values at the test inputs. For fixed hyperparameters, the joint
distribution of the training observations and test function values is
\[
\begin{bmatrix}
\mathbf y_{\mathrm{train}}\\
\mathbf f_{\mathrm{test}}
\end{bmatrix}
\sim
\mathcal N
\left(
\mathbf 0,
\begin{bmatrix}
\bK+\sI & \bK_*\\
\bK_*^T & \bK_{**}
\end{bmatrix}
\right).
\]
Conditioning on the training observations gives
\begin{align*}
p(
\mathbf f_{\mathrm{test}}
\mid
\mathbf X_{\mathrm{test}},
\mathbf X_{\mathrm{train}},
\mathbf y_{\mathrm{train}}
)
&=
\mathcal N(
\mathbf f_{\mathrm{test}}
\mid
\boldsymbol{\mu}_*,
\bm{\Sigma}_*
),\\
\boldsymbol{\mu}_*
&=
\bK_*^T
(\bK+\sI)^\inv
\mathbf y_{\mathrm{train}},\\
\bm{\Sigma}_*
&=
\bK_{**}
-
\bK_*^T
(\bK+\sI)^\inv
\bK_*.
\end{align*}
The posterior mean $\boldsymbol{\mu}_*$ predicts the latent function
values, while $\bm{\Sigma}_*$ quantifies their uncertainty conditional
on the chosen hyperparameters. For noisy test observations with
independent Gaussian noise of variance $\sigma^2$, the predictive
covariance is
$\bm{\Sigma}_*+\sigma^2\mathbf I_{n_*}$.
Fig.~\ref{fig:gp} illustrates the prior and posterior of a
one-dimensional GP.

In the following sections, we also use the shorthand
$\bX=\mathbf X_{\mathrm{train}}$,
$\by=\mathbf y_{\mathrm{train}}$, and
$\bX_*=\mathbf X_{\mathrm{test}}$.

\subsection*{Kernel model}

The kernel determines the covariance structure of the GP. We begin
with the squared-exponential kernel
\[
k(x,x')
=
\sigma_f^2
\exp\left(
-\dfrac{\|x-x'\|_2^2}{2\ell^2}
\right),
\]
where $\ell>0$ is the length scale and $\sigma_f>0$ is the signal
scale, so that $\sigma_f^2$ is the signal variance. This kernel is
symmetric and positive definite.

Define the squared Euclidean distance matrix
$\bS\in\R^{n\times n}$ by
\[
S_{pq}
=
\|x_p-x_q\|_2^2,
\qquad
p,q=1,\dots,n.
\]
The corresponding kernel matrix is
\[
\bK(\ell,\sigma_f)
=
\sigma_f^2
\exp\left(
-\dfrac{\bS}{2\ell^2}
\right),
\]
where the exponential is applied entrywise.

In our experiments, we use the squared-exponential kernel with
\emph{automatic relevance determination} (ARD), which assigns a
separate positive length scale $\ell_j$ to each input coordinate:
\begin{equation}
\label{eq:ard-kernel}
k(x,x')
=
\sigma_f^2
\exp\left(
-\dfrac{1}{2}
\sum_{j=1}^d
\dfrac{(x_j-x_j')^2}{\ell_j^2}
\right).
\end{equation}
The hyperparameter vector is
\[
\theta
=
[\ell_1,\dots,\ell_d,\sigma_f,\sigma]^T
\]
and contains $d+2$ parameters. The noise scale $\sigma$ enters the
observation covariance $\bK+\sI$. A large $\ell_j$ relative to the variation in the $j$th
input coordinate makes the kernel less sensitive to changes along
that coordinate.

For each coordinate $j=1,\dots,d$, define
$\bS^{(j)}\in\R^{n\times n}$ by
\[
S_{pq}^{(j)}
=
(x_{p,j}-x_{q,j})^2,
\qquad
p,q=1,\dots,n,
\]
so that
\[
\bS
=
\sum_{j=1}^d \bS^{(j)}.
\]
The ARD kernel matrix can then be written as
\[
\bK(\theta)
=
\sigma_f^2
\exp\left(
-\dfrac{1}{2}
\sum_{j=1}^d
\dfrac{\bS^{(j)}}{\ell_j^2}
\right),
\]
with the exponential again applied entrywise. The isotropic kernel
is recovered by setting $\ell_j=\ell$ for all $j$.
The distance matrices provide convenient notation for the kernel and
its derivatives; the low-rank implementation described below forms
only the required blocks.

\subsection*{Large-scale computation}

The main computational difficulty in GP regression is that $\bK$ is
generally dense. For $\sigma>0$, the observation covariance
$\bK+\sI$ is symmetric positive definite in exact arithmetic.
A direct Cholesky factorization requires $\mathcal O(n^3)$ operations
and $\mathcal O(n^2)$ storage, which becomes prohibitive for large
datasets \cite{williams2006gaussian}. This has motivated approaches
based on inducing points, structured kernel approximations, and
matrix-free iterative methods.

Inducing-point methods use a smaller collection of inputs to construct
a low-rank representation of the covariance structure; see, e.g.,
\cite{quinonero2005unifying}. Examples include the
\emph{subset of regressors} (SOR) \cite{silverman1985some} and the
\emph{fully independent training conditional} (FITC) approximation
\cite{snelson2005sparse}. FITC augments the low-rank training
covariance with a diagonal correction that restores the marginal
variances at the training inputs. Structured kernel interpolation,
developed by Wilson and Nickisch \cite{wilson2015kernel}, provides
another way to approximate kernel matrices efficiently. Such
approximations can be combined with Krylov methods for linear solves
and with stochastic trace estimation for log-determinant computations
\cite{Dong}. These techniques are used in scalable GP software such
as GPyTorch \cite{gardner2018gpytorch}.

Efficient trace estimation is also important in these computations
\cite{Skor,ubaru2017fast,meyer2021hutch++,
cortinovis2021randomized,persson2022improved,
wenger2022preconditioning,xu2026}. Matrix-free approaches evaluate
the action of $\bK$ without explicitly storing the dense matrix
\cite{Wagner2027,higp}. For posterior-mean computations, the resulting
linear systems are closely related to those arising in kernel ridge
regression and can benefit from preconditioning, including low-rank
preconditioners; see, e.g.,
\cite{chen2025randomly,pgp,epperly2025embrace,xu2025}.

Our method uses \AnchorNet\ to select landmarks for a low-rank kernel
approximation and combines this approximation with the FITC diagonal
correction. The resulting training objective and its gradients are
evaluated using low-rank factors, diagonal quantities, and small dense
matrices, as developed in \Cref{sec:training}.

\section{Training Gaussian Processes}
\label{sec:training}

Exact GP training estimates the hyperparameters by minimizing the
negative log-marginal likelihood
\begin{equation}
\label{eq:exact-nll}
f(\theta)
=
\dfrac{1}{2}
\left(
n\log(2\pi)
+
\ld(\tKK)
+
\by^T\tKK^\inv\by
\right),
\qquad
\tKK=\bK+\sI.
\end{equation}
For large datasets, the main computational bottlenecks are the linear
solve with $\tKK$ required by the quadratic term and the evaluation
of $\ld(\tKK)$. A direct Cholesky factorization requires
$\mathcal O(n^3)$ operations and $\mathcal O(n^2)$ storage
\cite{williams2006gaussian,bishop2006pattern}, while alternative
large-scale approaches typically rely on iterative solves and
trace-estimation techniques
\cite{gardner2018gpytorch,dong2017scalable}.

To reduce these costs, we construct a low-rank approximation using
landmarks selected by \AnchorNet\ \cite{cai2022fast}. The corresponding
kernel blocks satisfy
\[
\bKKnr\in\R^{n\times r},
\qquad
\bKKrr\in\R^{r\times r},
\qquad
\bKrn=\bKKnr^T\in\R^{r\times n}.
\]
The construction of the landmarks and kernel blocks is described in
\Cref{sec::anchornet}. The landmark locations are distinct and remain
fixed throughout hyperparameter optimization.

In the single-precision implementation, the landmark kernel matrix
is stabilized before factorization:
\[
\overline{\bK}_{11}
=
\bKKrr+\epsilon\sigma_f^2\bI_r,
\]
where $\epsilon\geq0$ is selected as described in \Cref{sec:alg}.
For notational simplicity, in the remainder of this section,
$\bKKrr$ and $\bKKrrinv$ denote the stabilized matrix and its inverse,
respectively. With this convention, the low-rank approximation is
\begin{equation}
\label{eq:anchornet-approx}
\bK
\approx
\hat{\bK}
=
\bKKnr\bKKrrinv\bKrn.
\end{equation}
All derivatives below are taken with the landmark set and the
accepted value of $\epsilon$ held fixed. Changes in $\epsilon$
caused by the stabilization procedure are not differentiated.

\subsection{The FITC-corrected low-rank model}

The low-rank approximation can underestimate the marginal variances
at the training inputs. For nonnegative stabilization, it satisfies
\[
0\preceq\hat{\bK}\preceq\bK
\]
in exact arithmetic, and hence
\[
[\hat{\bK}]_{ii}
\leq
k(x_i,x_i).
\]
When $\epsilon=0$, the standard Nystr\"om approximation gives equality
at the landmark locations. Positive stabilization generally removes
this interpolation property and increases the diagonal deficit.

We use the fully independent training conditional (FITC) correction
\cite{snelson2005sparse} to restore the marginal variances. Let
$\boldsymbol{\Lambda}$ be the diagonal matrix whose entries are
\[
\Lambda_{ii}
=
k(x_i,x_i)-[\hat{\bK}]_{ii}.
\]
The approximate covariance of the noisy training observations is
\begin{equation}
\label{eq:fitc}
\bK_{\mathrm{FITC}}
=
\hat{\bK}+\bD,
\qquad
\bD
=
\boldsymbol{\Lambda}+\sI.
\end{equation}
Thus, $\hat{\bK}+\boldsymbol{\Lambda}$ approximates the latent
covariance $\bK$, while $\bK_{\mathrm{FITC}}$ approximates
$\tKK=\bK+\sI$. Since $\Lambda_{ii}\geq0$ and $\sigma>0$,
\[
D_{ii}\geq\sigma^2>0,
\]
so both $\bD$ and $\bK_{\mathrm{FITC}}$ are positive definite in
exact arithmetic.

The approximate negative log-marginal likelihood used for training is
\begin{equation}
\label{eq:fitc-nll}
\hat f(\theta)
=
\dfrac{1}{2}
\left(
n\log(2\pi)
+
\ld(\bK_{\mathrm{FITC}})
+
\by^T\bK_{\mathrm{FITC}}^\inv\by
\right).
\end{equation}
The uncorrected low-rank model is recovered by setting
$\boldsymbol{\Lambda}=0$.

\subsection{Efficient evaluation of the objective}

The inverse and log-determinant in \eqref{eq:fitc-nll} can be
evaluated using low-rank factors, diagonal operations, and
$r\times r$ factorizations. Define
\begin{equation}
\label{eq:ksmall}
\bK_{\mathrm{small}}
\coloneqq
\bI_r
+
\bKKrrinv\bKrn\bD^\inv\bKKnr.
\end{equation}
Applying the Sherman--Morrison--Woodbury formula gives
\begin{align}
\boldsymbol{\alpha}
\coloneqq
\bK_{\mathrm{FITC}}^\inv\by
={}&
\bD^\inv\by
-
\bD^\inv\bKKnr
\bK_{\mathrm{small}}^\inv
\bKKrrinv\bKrn
\bD^\inv\by.
\label{eq:fitc-alpha}
\end{align}
Similarly, the matrix determinant lemma gives
\begin{equation}
\label{eq:fitc-logdet}
\ld(\bK_{\mathrm{FITC}})
=
\sum_{m=1}^n \log D_{mm}
+
\ld(\bK_{\mathrm{small}}).
\end{equation}

Although $\bK_{\mathrm{small}}$ is generally nonsymmetric, it is
similar to a symmetric positive definite matrix. Let
\[
\bKKrr=\bL\bL^T,
\qquad
\bW=\bL^\inv\bKrn\in\R^{r\times n}.
\]
Then
\begin{equation}
\label{eq:symmetric-small}
\bM
=
\bL^T\bK_{\mathrm{small}}\bL^{-T}
=
\bI_r+\bW\bD^\inv\bW^T.
\end{equation}
Since $\bD$ is positive definite,
$\bM\succeq\bI_r$ in exact arithmetic. The similarity transformation
also gives
\[
\ld(\bK_{\mathrm{small}})
=
\ld(\bM).
\]
The inverse action in \eqref{eq:fitc-alpha} can equivalently be written
as
\[
\boldsymbol{\alpha}
=
\bD^\inv\by
-
\bD^\inv\bW^T
\bM^\inv
\bW\bD^\inv\by.
\]
Thus, all dense factorizations required for evaluating the objective
involve only $r\times r$ matrices. The formulas can be evaluated using
diagonal scaling and triangular solves, without explicitly forming
the inverse matrices.

\subsection{Gradient evaluation}

To optimize the hyperparameters, we differentiate
\eqref{eq:fitc-nll}. For any hyperparameter $\theta_i$, the standard
identities for the derivative of a log-determinant and a matrix
inverse give
\begin{equation}
\label{eq:fitc-grad}
\pdv{\hat f}{\theta_i}
=
\dfrac{1}{2}
\left[
\tr\left(
\bK_{\mathrm{FITC}}^\inv
\pdv{\bK_{\mathrm{FITC}}}{\theta_i}
\right)
-
\boldsymbol{\alpha}^T
\pdv{\bK_{\mathrm{FITC}}}{\theta_i}
\boldsymbol{\alpha}
\right].
\end{equation}
Using
\[
\bK_{\mathrm{FITC}}
=
\hat{\bK}+\bD,
\]
we separate the gradient into low-rank and diagonal contributions:
\begin{align}
\pdv{\hat f}{\theta_i}
={}&
\dfrac{1}{2}
\left[
\tr\left(
\bK_{\mathrm{FITC}}^\inv
\pdv{\hat{\bK}}{\theta_i}
\right)
-
\boldsymbol{\alpha}^T
\pdv{\hat{\bK}}{\theta_i}
\boldsymbol{\alpha}
\right]
\nonumber+
\dfrac{1}{2}
\sum_{m=1}^n
q_m
\pdv{D_{mm}}{\theta_i},
\label{eq:fitc-grad-split}
\end{align}
where
\begin{equation}
\label{eq:q-def}
q_m
=
[\bK_{\mathrm{FITC}}^\inv]_{mm}
-
\alpha_m^2.
\end{equation}

The diagonal of $\bK_{\mathrm{FITC}}^\inv$ can be obtained without
forming the full inverse. Using $\bW=\bL^\inv\bKrn$ and
\eqref{eq:symmetric-small},
\begin{equation}
\label{eq:fitc-inv-diag}
\operatorname{diag}(\bK_{\mathrm{FITC}}^\inv)
=
\operatorname{diag}(\bD^\inv)
-
\operatorname{diag}
\left(
\bD^\inv\bW^T\bM^\inv\bW\bD^\inv
\right).
\end{equation}
The second diagonal is evaluated from blocks of the low-rank factors,
without forming an $n\times n$ matrix. It remains to evaluate the
derivatives of $\hat{\bK}$ and the associated reduced quantities.

\paragraph{Length-scale derivatives}

For an ARD length scale $\ell_j$, all three factors in
\[
\hat{\bK}
=
\bKKnr\bKKrrinv\bKrn
\]
depend on $\ell_j$. The matrix product rule
\cite{giles2008collected} therefore gives
\begin{align}
\pdv{\hat{\bK}}{\ell_j}
={}&
\pdv{\bKKnr}{\ell_j}
\bKKrrinv\bKrn
+
\bKKnr
\pdv{\bKKrrinv}{\ell_j}
\bKrn
+
\bKKnr\bKKrrinv
\pdv{\bKrn}{\ell_j}.
\label{eq:anchornet-product-deriv}
\end{align}
Using the ARD kernel \eqref{eq:ard-kernel} and the componentwise
squared-distance matrices $\bS^{(j)}$ defined in \Cref{sec:gp}, the
factor derivatives are
\begin{align}
\pdv{\bKKnr}{\ell_j}
&=
\dfrac{1}{\ell_j^3}
\bS^{(j)}(:,\mathrm{nyst})
\had
\bKKnr,
\nonumber\\
\pdv{\bKKrrinv}{\ell_j}
&=
-\dfrac{1}{\ell_j^3}
\bKKrrinv
\left(
\bS^{(j)}(\mathrm{nyst},\mathrm{nyst})
\had
\bKKrr
\right)
\bKKrrinv,
\nonumber\\
\pdv{\bKrn}{\ell_j}
&=
\dfrac{1}{\ell_j^3}
\bS^{(j)}(\mathrm{nyst},:)
\had
\bKrn,
\qquad
j=1,\dots,d.
\label{eq:ard-deriv}
\end{align}
Here $\had$ denotes the Hadamard product and $\mathrm{nyst}$ denotes
the landmark indices identified by \AnchorNet. With $\epsilon$ fixed,
the stabilization term $\epsilon\sigma_f^2\bI_r$ is independent of
$\ell_j$. Moreover, the diagonal of
$\bS^{(j)}(\mathrm{nyst},\mathrm{nyst})$ is zero, so its Hadamard
product with the stabilization term vanishes. Thus,
\eqref{eq:ard-deriv} remains valid under our convention that
$\bKKrr$ denotes the stabilized matrix.

The derivative in \eqref{eq:anchornet-product-deriv} admits the
low-rank factorization
\begin{equation}
\label{eq:uvw-ard}
\pdv{\hat{\bK}}{\ell_j}
=
\underbrace{
\begin{bmatrix}
\pdv{\bKKnr}{\ell_j}
&
\bKKnr
&
\bKKnr
\end{bmatrix}
}_{\mathbf U^{(j)}}
\underbrace{
\begin{bmatrix}
\bKKrrinv & & \\
& \pdv{\bKKrrinv}{\ell_j} & \\
& & \bKKrrinv
\end{bmatrix}
}_{\mathbf G^{(j)}}
\underbrace{
\begin{bmatrix}
\bKrn\\
\bKrn\\
\pdv{\bKrn}{\ell_j}
\end{bmatrix}
}_{\mathbf V^{(j)}}.
\end{equation}
The factors have dimensions
\[
\mathbf U^{(j)}\in\R^{n\times3r},
\qquad
\mathbf G^{(j)}\in\R^{3r\times3r},
\qquad
\mathbf V^{(j)}\in\R^{3r\times n}.
\]
This block representation can be evaluated termwise without
explicitly assembling the three factors. In particular,
\begin{equation}
\label{eq:quad-ard}
\boldsymbol{\alpha}^T
\pdv{\hat{\bK}}{\ell_j}
\boldsymbol{\alpha}
=
\boldsymbol{\alpha}^T
\mathbf U^{(j)}
\mathbf G^{(j)}
\mathbf V^{(j)}
\boldsymbol{\alpha},
\end{equation}
so the quadratic-form contribution requires no $n\times n$
derivative matrix.

The low-rank trace contribution can also be reduced to an
$r\times r$ computation. Since the derivative of $\bD$ is accounted
for separately in \eqref{eq:fitc-grad-split}, we differentiate
$\bK_{\mathrm{small}}$ with $\bD$ held fixed:
\begin{align}
\dK_{\mathrm{small,LR}}^{(\ell_j)}
\coloneqq{}&
\left.
\pdv{\bK_{\mathrm{small}}}{\ell_j}
\right|_{\bD\ {\rm fixed}}
\nonumber\\
={}&
\pdv{\bKKrrinv}{\ell_j}
\bKrn\bD^\inv\bKKnr
+
\bKKrrinv
\pdv{\bKrn}{\ell_j}
\bD^\inv\bKKnr
\nonumber+
\bKKrrinv\bKrn\bD^\inv
\pdv{\bKKnr}{\ell_j}.
\label{eq:dksmall-ard}
\end{align}
Differentiating the determinant identity with $\bD$ fixed gives
\begin{equation}
\label{eq:trace-ard}
\tr\left(
\bK_{\mathrm{FITC}}^\inv
\pdv{\hat{\bK}}{\ell_j}
\right)
=
\tr\left(
\bK_{\mathrm{small}}^\inv
\dK_{\mathrm{small,LR}}^{(\ell_j)}
\right).
\end{equation}

The diagonal contribution follows from
$k(x_m,x_m)=\sigma_f^2$, which is independent of the length scales:
\begin{equation}
\label{eq:dD-ell}
\pdv{D_{mm}}{\ell_j}
=
-
\left[
\pdv{\hat{\bK}}{\ell_j}
\right]_{mm}.
\end{equation}
Combining these expressions yields
\begin{align}
\pdv{\hat f}{\ell_j}
=
\dfrac{1}{2}
\Bigg[
&
\tr\left(
\bK_{\mathrm{small}}^\inv
\dK_{\mathrm{small,LR}}^{(\ell_j)}
\right)
-
\boldsymbol{\alpha}^T
\mathbf U^{(j)}
\mathbf G^{(j)}
\mathbf V^{(j)}
\boldsymbol{\alpha}
\nonumber\\
&
-
\sum_{m=1}^n
q_m
\left[
\mathbf U^{(j)}
\mathbf G^{(j)}
\mathbf V^{(j)}
\right]_{mm}
\Bigg],
\qquad
j=1,\dots,d.
\label{eq:ard-grad}
\end{align}

\paragraph{Signal-scale derivative}

The signal variance $\sigma_f^2$ scales every kernel block. Because
the stabilization has the form $\epsilon\sigma_f^2\bI_r$ with fixed
$\epsilon$, the stabilized landmark matrix also scales
proportionally to $\sigma_f^2$. Consequently,
\begin{equation}
\label{eq:sigmaf-factor-deriv}
\pdv{\bKKnr}{\sigma_f}
=
\dfrac{2}{\sigma_f}\bKKnr,
\qquad
\pdv{\bKrn}{\sigma_f}
=
\dfrac{2}{\sigma_f}\bKrn,
\qquad
\pdv{\bKKrrinv}{\sigma_f}
=
-\dfrac{2}{\sigma_f}\bKKrrinv.
\end{equation}
The product rule then gives
\begin{equation}
\label{eq:sigmaf-khat}
\pdv{\hat{\bK}}{\sigma_f}
=
\dfrac{2}{\sigma_f}\hat{\bK}.
\end{equation}

As in the length-scale calculation, the low-rank contribution holds
$\bD$ fixed:
\begin{equation}
\label{eq:dksmall-sigmaf}
\dK_{\mathrm{small,LR}}^{(\sigma_f)}
\coloneqq
\left.
\pdv{\bK_{\mathrm{small}}}{\sigma_f}
\right|_{\bD\ {\rm fixed}}
=
\dfrac{2}{\sigma_f}
\bKKrrinv\bKrn\bD^\inv\bKKnr
=
\dfrac{2}{\sigma_f}
\left(
\bK_{\mathrm{small}}-\bI_r
\right).
\end{equation}
Since
\[
\Lambda_{mm}
=
\sigma_f^2-[\hat{\bK}]_{mm},
\]
the diagonal derivative is
\begin{equation}
\label{eq:dD-sigmaf}
\pdv{D_{mm}}{\sigma_f}
=
\dfrac{2}{\sigma_f}\Lambda_{mm}.
\end{equation}
The resulting signal-scale gradient is
\begin{align}
\pdv{\hat f}{\sigma_f}
=
\dfrac{1}{2}
\Bigg[
&
\tr\left(
\bK_{\mathrm{small}}^\inv
\dK_{\mathrm{small,LR}}^{(\sigma_f)}
\right)
-
\dfrac{2}{\sigma_f}
\boldsymbol{\alpha}^T
\hat{\bK}
\boldsymbol{\alpha}
\nonumber
+
\dfrac{2}{\sigma_f}
\sum_{m=1}^n
q_m\Lambda_{mm}
\Bigg].
\label{eq:sigmaf-grad}
\end{align}

\paragraph{Noise derivative}

The low-rank approximation $\hat{\bK}$ and the diagonal correction
$\boldsymbol{\Lambda}$ are independent of the noise scale $\sigma$,
while
\[
\pdv{\bD}{\sigma}
=
2\sigma\bI_n.
\]
Thus, only the diagonal contribution in
\eqref{eq:fitc-grad-split} remains:
\begin{equation}
\label{eq:sigma-grad}
\pdv{\hat f}{\sigma}
=
\sigma\sum_{m=1}^n q_m
=
\sigma
\left(
\tr(\bK_{\mathrm{FITC}}^\inv)
-
\boldsymbol{\alpha}^T\boldsymbol{\alpha}
\right).
\end{equation}

All quantities in
\eqref{eq:ard-grad}, \eqref{eq:sigmaf-grad}, and
\eqref{eq:sigma-grad} can be evaluated using low-rank factors,
diagonal quantities, and reduced $r\times r$ matrices. Neither the
full kernel matrix nor any $n\times n$ derivative matrix needs to be
formed. The componentwise squared-distance matrices $\bS^{(j)}$
are likewise not stored explicitly; the required Hadamard products
are generated blockwise together with the corresponding kernel
blocks.

\subsection{Positive parameterization}

To enforce positivity without constrained optimization, we introduce
unconstrained parameters
\[
\boldsymbol{\zeta}
=
(\zeta_{\ell_1},\ldots,\zeta_{\ell_d},
 \zeta_{\sigma_f},\zeta_\sigma)^T
\in\R^{d+2}
\]
and parameterize
\[
\ell_j
=
\operatorname{softplus}(\zeta_{\ell_j}),
\qquad
\sigma_f
=
\operatorname{softplus}(\zeta_{\sigma_f}),
\]
and
\begin{equation}
\label{eq:softplus}
\sigma
=
10^{-4}
+
\operatorname{softplus}(\zeta_\sigma),
\qquad
\operatorname{softplus}(\zeta)
=
\log(1+e^\zeta).
\end{equation}
Since
\[
\operatorname{softplus}'(\zeta)
=
(1+e^{-\zeta})^{-1},
\]
the derivatives with respect to the unconstrained optimization
variables follow from the chain rule. Writing $\zeta_i$ for the
component corresponding to $\theta_i$, we obtain
\begin{equation}
\label{eq:softplus-grad}
\pdv{\hat f}{\zeta_i}
=
\pdv{\hat f}{\theta_i}
(1+e^{-\zeta_i})^{-1},
\qquad
i=1,\dots,d+2.
\end{equation}
This formula also applies to the noise parameter because its
additional floor is constant.

\section{\AnchorNet\ Low-Rank Approximation}
\label{sec::anchornet}

We use \AnchorNet\ \cite{cai2022fast} to select landmark points for
the low-rank kernel approximation. Given a set of landmarks $S$,
the standard Nystr\"om approximation is
\[
\bK
\approx
\hat{\bK}
=
\bKKnr\bKKrrinv\bKrn,
\]
where
$\bKKnr\in\R^{n\times r}$,
$\bKKrr\in\R^{r\times r}$, and
$\bKrn=\bKKnr^T\in\R^{r\times n}$.
Here $r=\#S$ is the number of distinct landmarks.

Let
\[
\mathcal X
=
\{x_1,\dots,x_n\}
\subset\R^d
\]
denote the training locations stored in the rows of
$\mathbf X_{\mathrm{train}}$.
Let $m$ be a positive integer specifying the target net size.
The construction has two stages: it first generates anchor points
adapted to the data geometry and then maps them to training points.
All nearest-point ties are resolved by a fixed rule.

\emph{Anchor construction.}
\begin{enumerate}
  \item Compute the smallest axis-aligned bounding box $B_0$
  containing $\mathcal X$ and generate an initial low-discrepancy set
  \[
  \mathcal T
  =
  \{t_1,\dots,t_s\}
  \subset B_0,
  \qquad
  s=\mathcal O(m).
  \]
  Possible constructions include Halton points and the adaptive
  tensor grids described in \cite{cai2022fast}.

  \item Assign each training point $x_j$ to one nearest point in
  $\mathcal T$ with respect to the $\ell_\infty$ norm. Specifically,
  choose
  \[
  a(j)
  \in
  \argmin_{1\leq i\leq s}
  \|x_j-t_i\|_\infty
  \]
  according to the fixed tie-breaking rule and define
  \[
  G_i
  =
  \{x_j:a(j)=i\},
  \qquad
  i=1,\dots,s.
  \]
  Discard empty groups and relabel the remaining groups as
  $G_1,\dots,G_Q$.

  \item For each nonempty group $G_i$, compute its smallest
  axis-aligned bounding box $B_i$. Let $\lambda(B_i)$ denote its
  $d$-dimensional Lebesgue measure. For the volume-based allocation,
  assume
  \[
  \sum_{k=1}^Q\lambda(B_k)>0
  \]
  and set
  \[
  m_i
  =
  \left\lceil
  m\,
  \frac{\lambda(B_i)}
  {\sum_{k=1}^Q\lambda(B_k)}
  \right\rceil.
  \]
  When $m_i>0$, generate a low-discrepancy set
  $\mathcal A^{(i)}\subset B_i$ with $m_i$ points.
  If $m_i=0$, set $\mathcal A^{(i)}=\varnothing$.

  \item Form the anchor net
  \[
  \mathcal A_X
  =
  \bigcup_{i=1}^Q
  \mathcal A^{(i)}.
  \]
\end{enumerate}

The ceiling allocation gives
\[
\#\mathcal A_X
\leq
\sum_{i=1}^Q m_i
\leq
m+Q
=
\mathcal O(m),
\]
because $Q\leq s=\mathcal O(m)$. Thus, the number of distinct anchors
need not equal the target size $m$.

\emph{From anchors to landmarks.}
For each anchor $y\in\mathcal A_X$, select a nearest training point,
\[
s(y)
\in
\argmin_{x\in\mathcal X}
\|x-y\|_\infty,
\]
using the fixed tie-breaking rule. The landmark set is
\[
S
=
\{s(y):y\in\mathcal A_X\}
\subset\mathcal X,
\]
with repeated selections removed. Consequently,
\[
r
=
\#S
\leq
\#\mathcal A_X
\leq
m+Q
=
\mathcal O(m),
\]
and $r\leq n$.

The volume-based allocation is designed to distribute anchors across
the data geometry. Mapping these anchors to the data produces
landmarks intended to improve the approximation quality and
numerical behavior of the Nystr\"om approximation
\cite{cai2022fast}.

The landmark-selection stage costs $\mathcal O(mdn)$ and depends
only on the input locations, requiring neither kernel matrix entries
nor kernel matrix--vector products \cite{cai2022fast}.
The selected landmarks define the low-rank approximation described
in \Cref{sec:training}, with stabilization implemented as described
in \Cref{sec:alg}.

\section{Algorithm}
\label{sec:alg}

\Cref{alg:anchornet-gp} summarizes the training and prediction
procedure. The landmarks are selected once and remain fixed during
optimization, with $r=\#S$ denoting the resulting approximation rank.
Initialization and optimizer settings are given in
\Cref{sec::appendix}.

\begin{algorithm}[t]
\caption{FITC Gaussian process training with \AnchorNet\ landmarks}
\label{alg:anchornet-gp}
\begin{algorithmic}[1]
\Require Training data $(\bX,\by)$; test inputs $\bX_*$;
target net size $m$; iteration budget $T\geq2$;
initial learning rate $\eta_0$;
low-rank predictive-mean flag $b_{\mathrm{LR}}$
\Ensure Hyperparameters $\theta$ and predictive mean
$\boldsymbol{\mu}_*$

\State $S\gets\textproc{AnchorNet}(\bX,m)$; $r\gets\#S$
\State $\theta_0\gets\textproc{Initialize}(\bX,\by)$
\State $\boldsymbol{\zeta}\gets
       \textproc{InversePositiveMap}(\theta_0)$
\State initialize the Adam optimizer state

\For{$t=0,\dots,T-1$}
  \State set $\eta_t$ using \eqref{eq:optimizer}
  \State $\theta\gets
         \textproc{PositiveMap}(\boldsymbol{\zeta})$
  \State form $\bKKrr=k(S,S)$,
         $\bKKnr=k(\bX,S)$, and $\bKrn=\bKKnr^T$
  \State $(\bL,\epsilon)\gets
         \textproc{StabilizedCholesky}(\bKKrr,\sigma_f)$
  \State compute $\bW,\bD,\bM$ as in \Cref{sec:training}
         and factor $\bM$
  \State evaluate $\hat f$ and $\nabla_\theta\hat f$
         using the stabilized factors, with $\epsilon$ fixed
  \State compute $\mathbf g=\nabla_{\boldsymbol{\zeta}}\hat f$
         using \eqref{eq:softplus-grad}
  \State $\boldsymbol{\zeta}\gets
         \textproc{Adam}
         \bigl(\boldsymbol{\zeta},
         \operatorname{clip}_{10}(\mathbf g),\eta_t\bigr)$
\EndFor

\State $\theta\gets
       \textproc{PositiveMap}(\boldsymbol{\zeta})$
\If{$b_{\mathrm{LR}}=\mathrm{True}$}
  \State recompute $\bL,\bW,\bD,\bM$ at the final $\theta$
  \State factor $\bM$ and solve
         $\bM\mathbf u=\bW\bD^\inv\by$
  \State solve $\bL^T\mathbf c=\mathbf u$
  \State evaluate
         $\boldsymbol{\mu}_*\gets k(\bX_*,S)\mathbf c$
\Else
  \State solve
         $\bigl(k(\bX,\bX)+\sigma^2\bI_n\bigr)
         \boldsymbol{\alpha}_{\mathrm{full}}=\by$
         at the final $\theta$
  \State evaluate
         $\boldsymbol{\mu}_*
         \gets k(\bX_*,\bX)\boldsymbol{\alpha}_{\mathrm{full}}$
\EndIf
\end{algorithmic}
\end{algorithm}

Training is performed in fp32. The stabilization routine factors
\[
\overline{\bK}_{11}
=
\bKKrr+\epsilon\sigma_f^2\bI_r
=
\bL\bL^T,
\]
starting with $\epsilon=10^{-4}$ and increasing it by a factor of
ten after an unsuccessful attempt, up to seven attempts.
Failure is reported if no acceptable factor is obtained; this did
not occur in the reported experiments. The gradient formulas in
\Cref{sec:training} use $\overline{\bK}_{11}$ and its inverse,
with the accepted $\epsilon$ held fixed. Final factor reconstruction
repeats the same stabilization procedure. The operator
$\operatorname{clip}_{10}$ clips the Euclidean norm of the
unconstrained gradient to $10$.

The factor computations are blocked by default. We select the
unblocked implementation only when an initial timing test shows
that it is at least $20\%$ faster and sufficient memory is
available. Kernel products during prediction can also be evaluated
in blocks.

The low-rank predictive mean uses only reduced solves and kernel
products with the landmarks. For the full-kernel mean, datasets with
fewer than $60{,}000$ training points use a dense factorization;
larger datasets use preconditioned CG with KeOps matrix--vector
products. The default fp32 preconditioner uses the first
$\min(r,200)$ landmarks. Additional fp64 solves establish the
full-kernel reference predictions on 3DRoad, Buzz, and HouseElectric,
where fp32 does not attain sufficient accuracy; see
\Cref{sec:experiments_new,sec:solve}.

\section{Numerical experiments}
\label{sec:experiments_new}

We evaluate the low-rank model described in
\Cref{alg:anchornet-gp}, which we refer to as ANGP
(\AnchorNet\ Gaussian process),\footnote{The code and all results 
needed to reproduce the experiments are available at \url{https://github.com/moritzp23/AnchorNet-GP}.} on the nine UCI 
regression datasets used by Maddox et al.~\cite{maddox2022low}. 
The datasets are obtained
through the public \texttt{bayesian\_benchmarks} interface and range
from $n=13{,}500$ to $n=1{,}844{,}352$ training points, with input
dimension $d$ ranging from $3$ to $90$, as summarized in
\Cref{tab:summary}. We use the standard preprocessing provided by
\texttt{bayesian\_benchmarks}; in particular, the input features and
targets are standardized.

Throughout the experiments, we use the same kernel, hyperparameter
initialization, optimizer, and optimization budget for all datasets.
The only model parameter varied across datasets is the number of
\AnchorNet\ landmarks, which determines the approximation rank $r$.
We use the squared-exponential kernel with automatic relevance
determination (ARD) in \eqref{eq:ard-kernel}, as in
\cite{maddox2022low}, resulting in $d+2$ hyperparameters.

The hyperparameters are initialized using the data-driven procedure
described in \Cref{sec::appendix}: the length scales are initialized
using the median pairwise distance heuristic, the signal scale using
the standard deviation of the target, and the noise scale using the
mode of its prior. We optimize the FITC negative log-marginal
likelihood using Adam for a fixed budget of $100$ iterations with a
cosine-decayed learning rate. We use neither a validation split nor
early stopping and report the final iterate for every run. With the
initialization and learning-rate schedule used here, we did not
observe divergence, collapse of a length scale, or activation of the
noise floor in the reported experiments.

Each evaluation of the FITC objective and its gradient has cost
$\mathcal O(nr(r+d))$ in our implementation. The unblocked
implementation stores $\mathcal O(nr)$ intermediate quantities,
whereas the blocked implementation processes the training data in
row blocks of size $B$ and requires
$\mathcal O(rB+r^2)$ working memory. At the beginning of each fit, we
time one objective-and-gradient evaluation using both implementations.
The blocked implementation is used by default; the unblocked variant
is selected only when it is at least $20\%$ faster and fits in
available memory. Additional implementation details are given in
\Cref{sec::appendix}.

For prediction, datasets with fewer than $60{,}000$ training points
use a dense full-kernel factorization. For the larger datasets, the
full-kernel predictive mean is computed using preconditioned conjugate
gradients (CG), with matrix--vector products evaluated through KeOps
\cite{JMLR:v22:20-275}. We also evaluate the predictive mean obtained
directly from the fitted low-rank model, which avoids the full-kernel
linear solve.

Each configuration is evaluated on five data splits, and we report
the mean and standard deviation. For each dataset, we test several
landmark counts and report every tested configuration in
Table~\ref{tab:rmse-rank}. When we subsequently identify the cheapest
rank that reaches a reference RMSE, this should be interpreted as a
retrospective cost--accuracy comparison rather than as a rank-selection
procedure used during training.

All ANGP timings are measured on an NVIDIA A40 GPU with $48$\,GB of
memory. We compare against the V100 ($32$\,GB) timings reported by
Maddox et al.~\cite{maddox2022low}, denoted by \texttt{Maddox} in the
tables. Following their experimental setup, the reported fitting times
include hyperparameter optimization but exclude the predictive solve.

The two sets of timings are obtained on different hardware. The A40
has roughly a $2.4\times$ advantage in nominal single-precision
throughput over the V100, so part of the observed speedup is due to
hardware rather than the algorithm alone. We therefore report the raw
timings and interpret a speedup divided by $2.4$ as an approximate
hardware-normalized comparison. Under this normalization, ANGP remains
faster than the reference on all datasets, although the advantage on
3DRoad becomes small.

The remainder of this section examines three aspects of the method:
the cost--accuracy tradeoff as the approximation rank varies, the
difference between the low-rank and full-kernel predictive means, and
the numerical behavior of the full-kernel predictive solve in single
precision.

\subsection{Comparison with the reference}

Table~\ref{tab:summary} reports, for each dataset, the cheapest tested
rank whose mean test RMSE reaches the reference value reported in
\cite{maddox2022low}. As noted above, this selection is performed
retrospectively from the tested ranks and is intended to summarize the
cost required to attain the reference accuracy.

Whenever the fp32 full-kernel solve converges, we compute the reported
RMSE using that predictive mean. For 3DRoad, Buzz, and HouseElectric,
the fp32 CG solve does not attain sufficient accuracy because of
finite-precision effects. For these datasets, we instead report the
RMSE obtained from the corresponding fp64 full-kernel solve; these
entries are marked by $^{\ast}$ in \Cref{tab:summary}.

The Protein RMSE reported in \cite[Table~1]{maddox2022low} is $0.055$.
This value is inconsistent with both the half-precision result
($0.635$) in the same table and the values $0.481$--$0.509$ reported
in the appendix of that paper. We therefore interpret $0.055$ as a
misplaced decimal point and use $0.55$ as the reference value
throughout.

\input{tables/summary.tex}

Table~\ref{tab:rmse-rank} shows that ANGP reaches the reference RMSE
on all nine datasets for at least one of the tested ranks. Increasing
$r$ generally improves the quality of the low-rank model but also
increases the fitting cost, which grows quadratically with $r$ in the
dominant low-rank operations. The rank required to reach the reference
accuracy varies substantially across datasets.

For example, 3DRoad requires $r=3000$ to reach the reference accuracy.
At $r=6000$, ANGP is already slower than \texttt{Maddox} even before
accounting for the hardware difference; see
Table~\ref{tab:cost-3droad}. Across the datasets, the raw fitting
speedups range from $3.4\times$ on 3DRoad to $1120\times$ on Buzz.
On 3DRoad, the observed speedup is only somewhat larger than the 
approximately $2.4\times$ hardware advantage of the A40, leaving 
about $1.4\times$ after this approximate hardware normalization. 
One contributing factor is that 3DRoad requires a relatively large 
approximation rank while having only $d=3$ input dimensions, making 
the kernel matrix--vector products used by the reference method 
comparatively inexpensive.

\subsection{The effect of rank}

Table~\ref{tab:rmse-rank} reports the test RMSE for every rank tested.
Entries that match or improve upon the reference value in
\cite[Table~1]{maddox2022low} are shown in bold. The rank required to
reach the reference varies by roughly two orders of magnitude and
does not exhibit a simple dependence on either the number of training
points $n$ or the input dimension $d$.

As in the previous subsection, the full-kernel predictive mean is
computed in fp64 for 3DRoad, Buzz, and HouseElectric when the fp32 CG
solve does not attain sufficient accuracy. On Buzz, increasing the
model rank does not lead to a corresponding improvement in the
full-kernel predictive RMSE, and some of the best results occur at
relatively small ranks.

\input{tables/rmse_vs_rank}
\input{tables/rmse_vs_rank_lr}

The reference RMSE is reached at rank $100$ for PoleTele and Bike,
at rank $1600$ for Kin40K, and at rank $3000$ for both Song and
3DRoad. The variation cannot be explained by $n$ or $d$ alone. For
example, Song, with $n=464{,}000$, requires rank $3000$, whereas
HouseElectric, with $n=1.84$M, reaches the reference at rank $800$.
Similarly, 3DRoad has only $d=3$ input dimensions but requires a
larger rank than PoleTele, which has $d=26$. The required rank is
instead strongly influenced by how rapidly the fitted kernel spectrum
decays, together with the geometry of the particular dataset; these
properties are not known before fitting.

The rank also exposes a clear accuracy--cost tradeoff. On Bike, for
example, the RMSE decreases from $0.0707$ at rank $100$ to $0.0303$
at rank $1600$, while the fitting time increases from $2.4$\,s to
$10.2$\,s. Rank $100$ is sufficient to match the reference, but a
larger rank can substantially improve predictive accuracy at additional
cost. A similar effect occurs on Kin40K, where the RMSE decreases to
$0.0649$ at rank $8000$, compared with the reference value $0.099$,
although rank $1600$ is already sufficient to reach the reference.
Thus, the cheapest rank that matches the reference should not be
interpreted as an accuracy-optimal rank.

Table~\ref{tab:rmse-rank-lr} reports the same trained models but uses
the predictive mean of the low-rank FITC model instead of solving the
full-kernel system. This predictor requires no iterative full-kernel
solve, can be evaluated in closed form using the low-rank factors, and
therefore cannot suffer from iterative convergence failure. Across all
reported configurations, its RMSE has a median relative increase of
approximately $27\%$ compared with the full-kernel predictive mean.
It reaches the reference accuracy on Elevators, Bike, and Kin40K.
The gap generally decreases as $r$ increases, but at markedly
different rates across datasets. On HouseElectric, for example, the
low-rank mean is within approximately $4\%$ of the full-kernel result
at $r=800$, whereas on Buzz it remains at least $24\%$ worse over the
tested ranks.

Tables~\ref{tab:cost-3droad}--\ref{tab:cost-houseelectric} compare the
low-rank predictive mean with the fp32 and fp64 full-kernel predictive
solves on the three largest datasets for which these comparisons are
most informative. Results are averaged over five splits. The
``effective rank'' (eff.) denotes the numerical rank of $\bKKrr$;
when it is substantially smaller than the nominal landmark count,
additional landmarks contribute little independent numerical
information. The CG columns report predictive-solve times, whereas the
column $\times$ref reports the fitting-time speedup relative to
\texttt{Maddox}. Entries that reach the \texttt{Maddox} reference RMSE
are shown in bold.

\input{tables/cost_3droad}
\input{tables/cost_buzz}
\input{tables/cost_houseelectric}

The low-rank predictive mean requires at most $4.1$\,s in these
experiments, compared with $33$\,s to $25$\,min for the fp32
full-kernel CG solve and up to $24$\,h for the fp64 solve. Consequently,
for the low-rank predictor, fitting dominates prediction time, whereas
the full-kernel solve can itself become a dominant component of the
total cost.

HouseElectric illustrates the advantage of avoiding the full-kernel
solve. The fp32 CG solve does not attain the required accuracy at any
tested rank, while the low-rank mean reaches an RMSE of $0.0536$ at
$r=1000$, within the standard deviation of the reference
($0.052\pm0.002$). Its prediction requires $0.9$\,s following a
$544$\,s fit, and the fitting stage is $78.5\times$ faster than the
reported \texttt{Maddox} reference. The complete low-rank pipeline in
this case remains in single precision. On 3DRoad and Buzz, by
contrast, the low-rank predictive mean is inexpensive but does not
reach the accuracy of the full-kernel predictor, so the latter remains
necessary when the highest predictive accuracy is required.

\subsection{The predictive-mean solve in single precision}

For datasets with $n<60{,}000$, we compute the full-kernel predictive
mean using a dense factorization. For the four larger datasets, we
instead use preconditioned CG in fp32, with kernel matrix--vector
products evaluated using KeOps \cite{JMLR:v22:20-275}. The
preconditioner is constructed from the low-rank \AnchorNet\
approximation using at most $200$ landmarks.

The behavior of the fp32 solve varies substantially across datasets.
It converges on Song, stagnates at a relative residual of approximately
$10^{-2}$ on 3DRoad, varies strongly across splits on Buzz, and does
not attain the required accuracy on HouseElectric; see
Tables~\ref{tab:cost-3droad}--\ref{tab:cost-houseelectric}.

This behavior is closely related to conditioning. Since $\bK$ is
positive semidefinite,
\[
\kappa_2(\bK+\sigma^2\bI)
=
\frac{\lambda_{\max}(\bK)+\sigma^2}
     {\lambda_{\min}(\bK)+\sigma^2}.
\]
When $\lambda_{\min}(\bK)$ is small, reducing $\sigma^2$ can therefore
increase the condition number approximately in proportion to
$1/\sigma^2$. In finite precision, the attainable accuracy of an
iterative solve generally deteriorates as the product of the unit
roundoff and the condition number grows. Thus, the small fitted noise
levels are an important contributor to the observed loss of accuracy
in fp32. For example, the fitted values are approximately
$\sigma^2=6\times10^{-2}$ on Buzz and
$\sigma^2=2.5\times10^{-3}$ on HouseElectric.

The preconditioner can introduce an additional finite-precision
difficulty. Let
\[
r_p=\min(r,200)
\]
and let $\bW_p\in\R^{r_p\times n}$ denote the low-rank factor
constructed from the first $r_p$ landmarks. We compute a thin QR
factorization
\[
\begin{bmatrix}
\bW_p^T\\
\sigma\bI_{r_p}
\end{bmatrix}
=
\bQ\bR,
\qquad
\bQ_1=\bQ(1:n,:).
\]
The preconditioner is then applied in the form
\begin{equation}
P^{-1}\mathbf r
=
\sigma^{-2}
\left(
\mathbf r-\mathbf Q_1\mathbf Q_1^T\mathbf r
\right).
\end{equation}

When a substantial component of $\mathbf r$ lies in the span of
$\mathbf Q_1$, the two vectors
$\mathbf r$ and $\mathbf Q_1\mathbf Q_1^T\mathbf r$ can be close, so
their subtraction is susceptible to cancellation in single precision.
The resulting rounding error is subsequently multiplied by
$\sigma^{-2}$. This amplification factor is approximately $16$ on
Buzz and $400$ on HouseElectric. The observation is consistent with
\Cref{fig:he-precond}, where increasing the preconditioner rank on
HouseElectric can worsen the true residual.

Repeating the full-kernel solve in fp64 removes this accuracy barrier
in the cases considered here. However, the A40 has much lower fp64
than fp32 throughput, so these reference solves can require many hours.
We therefore use fp64 only when needed to establish the predictive
accuracy attainable by the full-kernel model. A more detailed analysis
of these solves is provided in \Cref{sec:solve}.

\section{Mixed-precision experiments}
\label{sec:outlook_new}

The training algorithm described above is implemented primarily in
fp32. Profiling shows that approximately two-thirds of the fitting
time is spent in dense matrix multiplications (GEMMs) and triangular
solves (TRSMs), making these operations natural candidates for
lower-precision arithmetic. We therefore investigate whether selected
operations can be performed with fp16 inputs while retaining fp32
accumulation. Our experiments show that the attainable speedup is
modest and, more importantly, that the resulting computation is not
reliable for all datasets.

In the implementation considered here, the triangular solves associated
with the Cholesky factors remain in fp32. Consequently, only four
large GEMMs in each objective-and-gradient evaluation are candidates
for reduced precision. One of these forms the matrix
\[
\bM
=
\bI_r+\bW\bD^{-1}\bW^T,
\]
which is similar to $\bK_{\mathrm{small}}$ and is the matrix factored
by Cholesky. In exact arithmetic,
\[
\lambda_{\min}(\bM)\geq 1.
\]
The remaining three GEMMs appear only in the gradient and therefore do
not directly affect either the objective value or the Cholesky
factorization.

We consider two mixed-precision variants, both using fp16 inputs with
fp32 accumulation. The first performs the GEMM forming $\bM$ and one
gradient GEMM in reduced precision; we refer to this as the
\emph{two-GEMM} variant. The second performs all four large GEMMs in
reduced precision. For the four-GEMM variant, the GEMM inputs are
scaled by powers of two to avoid overflow; without this scaling, one
gradient block exceeds the dynamic range of fp16 near the fitted
solution and produces \texttt{NaN} length-scale gradients. The
two-GEMM variant is useful for isolating the stability of the GEMM that
forms $\bM$, whereas the four-GEMM variant represents the largest
precision reduction considered without changing the triangular solves.

Table~\ref{tab:mp-fit-3droad} compares single fits on 3DRoad with
$r=3000$ using fp32 and the two mixed-precision variants. All runs use
split 73, the ARD kernel, and $100$ Adam iterations with the same
cosine learning-rate schedule. The two-GEMM variant produces nearly
the same fitted hyperparameters and test RMSE as the fp32 run:
$0.3261$ compared with $0.3258$, while reducing the fit time from
$465$\,s to $449$\,s, a speedup of only $1.04\times$. Using all four
GEMMs in reduced precision lowers the fit time to $424$\,s
($1.10\times$ speedup), with a test RMSE of $0.3277$.

Because these results are based on a single fit for each precision
configuration, the small differences in objective value and RMSE
should not be interpreted as systematic accuracy degradation. At the
fp32 fitted hyperparameters, however, we can compare the gradients
directly. The relative gradient error is $1.65\%$ for the two-GEMM
variant and $1.35\%$ for the four-GEMM variant. Thus, for this
particular problem, reducing the precision of these GEMMs introduces
only a modest perturbation to an individual gradient evaluation but
provides correspondingly modest performance gains.

\begin{table}[t]
  \centering
  \small
  \caption{Single fits on 3DRoad with $r=3000$ using fp32 and two
  mixed-precision variants.}
  \label{tab:mp-fit-3droad}
  \begin{adjustbox}{max width=\textwidth}
  \begin{tabular}{lrrr}
    \toprule
     & fp32 & 2 GEMMs in fp16 & 4 GEMMs in fp16 \\
    \midrule
    Fit time (s) & 465 & 449 & 424 \\
    Final train NLL & 0.2365 & 0.2373 & 0.2406 \\
    Median length scale & 0.03906 & 0.03906 & 0.03911 \\
    $\sigma_f$ & 0.5930 & 0.5907 & 0.5899 \\
    $\sigma$ & 0.2308 & 0.2308 & 0.2303 \\
    Effective rank & 2921 & 2919 & 2931 \\
    Test RMSE & 0.3258 & 0.3261 & 0.3277 \\
    Reduced-precision $\bM$ used & -- & 74 of 100 & 79 of 100 \\
    \bottomrule
  \end{tabular}
  \end{adjustbox}
\end{table}

The limited speedup is consistent with the operation mix of the
algorithm. The triangular solves, kernel evaluations, and several
other operations remain in fp32. For the two-GEMM variant, the speedup
of a single objective-and-gradient evaluation grows with rank, from no
measurable improvement at $r=250$ to approximately $1.13\times$ at
$r=6000$ on 3DRoad; see Table~\ref{tab:mp-scan}. At $r=3000$, using
all four GEMMs in reduced precision gives an evaluation speedup of
approximately $1.15\times$. At the level of a complete fit, however,
the observed speedup remains close to $1.1\times$.

\subsection{Stability of the reduced-precision factorization}

The main numerical difficulty is the GEMM used to form $\bM$. To
distinguish exact-arithmetic properties from finite-precision effects,
let
\[
\bM
=
\bI_r+\bW\bD^{-1}\bW^T
\]
denote the exact matrix,
let $\widetilde{\bM}$ denote the matrix obtained using the fp32
baseline computation, and let $\widehat{\bM}$ denote the corresponding
matrix when the selected GEMM is evaluated with fp16 inputs and fp32
accumulation.

Although $\lambda_{\min}(\bM)\geq1$ in exact arithmetic, rounding in
the matrix products can perturb the computed matrices. We write
\[
\widehat{\bM}
=
\bM+\Delta\bM.
\]
As a rough error model, the norm of $\Delta\bM$ can grow proportionally
to the fp16 unit roundoff and the scale of the matrix product. By
Weyl's inequality, positive definiteness of the computed matrix is
guaranteed whenever
\[
\|\Delta\bM\|_2
<
\lambda_{\min}(\bM).
\]
This suggests the heuristic condition
\[
u_{\mathrm{fp16}}\,\kappa_2(\bM)
\lesssim 1,
\]
although the actual behavior also depends on the structure and scaling
of the matrix multiplication.

Table~\ref{tab:mp-scan} reports one mixed-precision
objective-and-gradient evaluation per rank at hyperparameters obtained
from the fp32 fit on split 73. On 3DRoad, the reduced-precision
factorization succeeds at every tested rank between $50$ and $6000$,
even though
$u_{\mathrm{fp16}}\kappa_2(\bM)\approx12$ at $r=3000$.
Thus, the heuristic condition is conservative for this problem.

The behavior is substantially worse on Buzz and HouseElectric. Their
estimated condition numbers lie between approximately
$1.3\cdot10^6$ and $1.2\cdot10^7$. The reduced-precision
factorization fails for $7$ of the $8$ tested ranks on Buzz and for
$9$ of the $11$ tested ranks on HouseElectric. At the few ranks where
the factorization succeeds, the relative gradient errors are already
large. We therefore do not attempt complete reduced-precision fits on
these two datasets.

\begin{table}[t]
  \centering
  \small
  \caption{Single objective-and-gradient evaluations with the
  two-GEMM mixed-precision variant. The condition number is evaluated
  from the fp32 baseline matrix.}
  \label{tab:mp-scan}
  \begin{adjustbox}{max width=\textwidth}
  \begin{tabular}{lccccc}
    \toprule
    Dataset & $r$ & $\kappa_2(\widetilde{\bM})$
    & Factorization & Rel.\ gradient error
    & Time fp32 / mixed (s) \\
    \midrule
    3DRoad & 250  & $1.3\cdot10^{4}$ & succeeds & 0.029 & 0.16 / 0.17 \\
    3DRoad & 1000 & $2.1\cdot10^{4}$ & succeeds & 0.034 & 0.95 / 0.91 \\
    3DRoad & 3000 & $2.5\cdot10^{4}$ & succeeds & 0.016 & 4.58 / 4.20 \\
    3DRoad & 6000 & $2.5\cdot10^{4}$ & succeeds & 0.012 & 13.6 / 12.0 \\
    Buzz & 50 to 1000
         & $2.4\cdot10^{6}$ to $5.4\cdot10^{6}$
         & fails (7/8) & 4.2 & -- \\
    HouseElectric & 50 to 3000
         & $1.3\cdot10^{6}$ to $1.2\cdot10^{7}$
         & fails (9/11) & 3.0 and 29.6 & -- \\
    \bottomrule
  \end{tabular}
  \end{adjustbox}
\end{table}

These experiments also show that conditioning at the final fitted
hyperparameters does not characterize the entire optimization
trajectory. Table~\ref{tab:mp-trajectory} follows the 3DRoad
two-GEMM fit through selected Adam iterations. The exact-arithmetic
matrix $\bM$ remains positive definite by construction, but the
computed matrices can lose this property because of rounding.

Between iterations $18$ and $40$, the signal-to-noise ratio
$\sigma_f^2/\sigma^2$ ranges from approximately $25$ to $120$,
compared with $6.6$ near the final iterate. Over this part of the
trajectory, the largest eigenvalue of the fp32 baseline matrix grows
to between $2.2\cdot10^5$ and $1.2\cdot10^6$, while the
reduced-precision matrix $\widehat{\bM}$ becomes indefinite, with its
smallest eigenvalue between $-4.2$ and $-0.15$. At iterations $8$--$10$,
the reduced-precision GEMM instead overflows.

Because the reduced-precision matrix still admits a Cholesky
factorization at iteration $7$, despite a large spectral range, the
condition number is best interpreted as an indicator of numerical
risk rather than a deterministic breakdown criterion.

\begin{table}[t]
  \centering
  \small
  \caption{Computed spectra along the 3DRoad mixed-precision fit.
  $\widetilde{\bM}$ denotes the fp32 baseline matrix and
  $\widehat{\bM}$ the matrix obtained from the reduced-precision
  GEMM.}
  \label{tab:mp-trajectory}
  \begin{adjustbox}{max width=\textwidth}
  \begin{tabular}{lcccccc}
    \toprule
    Iteration
    & $\sigma_f$
    & $\sigma$
    & $\lambda_{\max}(\widetilde{\bM})$
    & $\lambda_{\min}(\widetilde{\bM})$
    & $\lambda_{\min}(\widehat{\bM})$
    & Factorization \\
    \midrule
    0
    & 1.00 & 1.414
    & $1.2\cdot10^{5}$ & 1.000 & 0.997 & succeeds \\
    7
    & 2.07 & 0.611
    & $1.5\cdot10^{6}$ & 1.000 & 0.975 & succeeds \\
    8 to 10
    & 2.25 to 2.61 & 0.55 to 0.57
    & $1.7\cdot10^{6}$ to $1.8\cdot10^{6}$
    & 0.999 & \texttt{Inf} & overflow \\
    16
    & 3.64 & 0.634
    & $5.9\cdot10^{5}$ & 0.936 & 0.222 & succeeds \\
    18 to 40
    & 4.07 to 1.44 & 0.60 to 0.30
    & $2.2\cdot10^{5}$ to $1.2\cdot10^{6}$
    & 0.87 to 0.99 & $-4.2$ to $-0.15$ & fails \\
    44
    & 1.01 & 0.291
    & $8.6\cdot10^{4}$ & 0.997 & 0.535 & succeeds \\
    Final
    & 0.59 & 0.231
    & $2.5\cdot10^{4}$ & 1.00 & -- & succeeds \\
    \bottomrule
  \end{tabular}
  \end{adjustbox}
\end{table}

To allow the mixed-precision fit to complete, we recompute the GEMM
forming $\bM$ in fp32 whenever the Cholesky factorization of
$\widehat{\bM}$ fails. Without this fallback, the 3DRoad fit terminates
at iteration $8$. With the fallback, the reduced-precision matrix is
accepted in $74$ of the $100$ Adam iterations for the two-GEMM variant
and in $79$ of $100$ iterations for the four-GEMM variant. The
power-of-two scaling used in the latter also prevents the observed
overflow. The fallback, however, detects only factorization failure;
it does not guarantee that a successful reduced-precision
factorization produces sufficiently accurate gradients. Instead of 
recomputing $\widehat{\bM}$ in fp32, the failed factorization can be 
repaired by raising the eigenvalues of $\widehat{\bM}$ below $1$ to 
$1$, which uses $\lambda_{\min}(\bM)\geq 1$. On 3DRoad, this keeps all 
four GEMMs in fp16 in every iteration and gives a test RMSE of $0.3259$,
compared with $0.3258$ in fp32. Such a modified factorization is 
therefore a possible fix when the loss of positive definiteness is the 
primary failure, but not on Buzz and HouseElectric, where the 
relative gradient error remains between $5$ and $650$ after the repair.

We also tested TF32, which retains a larger dynamic range but has
reduced mantissa precision, and observed behavior qualitatively
similar to fp16. This indicates that overflow is not the only
limitation: even when scaling prevents overflow, insufficient
precision can still perturb $\bM$ enough to compromise its
factorization or substantially degrade the gradient.

Buzz and HouseElectric are also the datasets for which $\bKKrr$ is
already close to numerically rank deficient in fp32, with observed
numerical ranks of $443$ out of $1000$ and $1724$ out of $3000$,
respectively. This is consistent with their greater sensitivity to
further precision reduction.

Overall, these experiments show that selected GEMMs in the
\AnchorNet--FITC training procedure can be executed below fp32 for
moderately conditioned problems, but the resulting performance gain is
small and the approach is not robust for ill-conditioned cases. In
particular, a Cholesky-based training method requires the computed
reduced matrix to remain positive definite, which cannot in general be
guaranteed after aggressive precision reduction. This observation is
consistent with the use of iterative methods rather than Cholesky
factorization for half-precision GP training in
\cite{maddox2022low}.

\section{Conclusion}

We developed a scalable Gaussian process training approach that combines
\AnchorNet\ landmark selection with an FITC correction and evaluates the
resulting objective and analytic gradients using low-rank factors and
small dense matrices. Across nine regression datasets, using the same
kernel, initialization, optimizer, and optimization budget while varying
the approximation rank, ANGP reaches the reference predictive accuracy
on all nine datasets. The corresponding fitting speedups range from
$3.4\times$ to $1120\times$, although these timings are obtained on an
A40 GPU with roughly a $2.4\times$ single-precision throughput advantage
over the V100 used for the reference results.

Our experiments also highlight an important distinction between training
and prediction in reduced precision. The low-rank training procedure can
be carried out effectively in single precision, whereas the full-kernel
predictive solve becomes numerically challenging on 3DRoad, Buzz, and
HouseElectric. On these datasets, the fp32 iterative solve does not
attain sufficient accuracy and an fp64 solve is required to establish
the full-kernel reference prediction. The low-rank predictive mean
avoids this solve entirely and requires only a few seconds. It reaches
the reference accuracy on Elevators, Bike, and Kin40K, and on
HouseElectric it reaches the reference within its reported standard
deviation while keeping the complete training-and-prediction pipeline in
single precision.

We further investigated whether selected dense matrix operations in the
training procedure can be reduced below fp32. Using fp16 or TF32 for
selected GEMMs provides only modest additional speedups and can become
unreliable when the reduced system is ill-conditioned. In particular,
rounding can destroy the positive definiteness required by the Cholesky
factorization even though the corresponding matrix is positive definite
in exact arithmetic. These results suggest that fp32 provides a useful
balance between performance and numerical robustness for the proposed
training procedure, while more aggressive precision reduction requires
additional stabilization.

The low-rank predictive mean is therefore an attractive first choice
when prediction cost is important. When its accuracy is insufficient,
the full-kernel predictive mean can be computed as a second stage, using
higher precision when required by the conditioning of the problem. This
suggests a practical adaptive strategy in which the low-rank prediction
is computed first and its accuracy is assessed on validation data before
deciding whether the more expensive full-kernel solve is necessary.
Developing a reliable criterion for this decision, together with
automatic selection of the approximation rank, is an interesting
direction for future work.
\section*{Acknowledgments}
The authors used Claude to assist with coding, including parameter
initialization, parameter sweeps, and plotting and table formatting,
and ChatGPT to assist with manuscript editing and revision. All
AI-assisted code and manuscript content were checked and verified by
the authors. The authors assume responsibility for all content.

% AFAIK we have to disclose the use of AI as specified in https://epubs.siam.org/artificial-intelligence, under a acknowledgments or declarations sections.
\bibliographystyle{siamplain}
\bibliography{references}	

\appendix

\section{Experimental setup}
\label{sec::appendix}

\subsection*{Initialization}
All datasets use the same initialization. The length scales start at the median pairwise distance of a random subsample $\mathcal{S}$ of $2{,}000$ training inputs,
\begin{equation}
  \ell_0 = \operatorname{median}\bigl\{\lVert x_i - x_j \rVert_2 : i < j,\ x_i, x_j \in \mathcal{S}\bigr\},
  \label{eq:median-heuristic}
\end{equation}
with the same value for all $d$ length scales. The signal scale starts at $\sigma_{f,0} = \operatorname{std}(y) = 1$, since the targets are standardized, and the noise at the mode of the $\mathrm{Gamma}(1.1, 0.05)$ prior on $\sigma^2$, i.e.\ $\sigma_0 = \sqrt{2} \approx 1.414$. A large initial noise keeps $\bK_{\mathrm{small}}$ well conditioned at the start; a small one multiplies the rounding errors by $\sigma^{-2}$ and was the main cause of failed fits in single precision. Each split is a single run from this starting point, without restarts.

\subsection*{Optimization}
We run Adam for $T = 100$ iterations with gradients clipped to norm $10$ and the learning rate
\begin{equation} \label{eq:optimizer}
  \eta_t = \eta_0\Bigl(0.02 + 0.98 \cdot \tfrac{1}{2}\bigl(1 + \cos(\pi t / (T-1))\bigr)\Bigr) \quad \text{and} \quad \eta_0 = 0.2 .
\end{equation}
The decay allows one learning rate for all datasets: 3DRoad has to move $4.1$ in $\log \ell$ from its starting point, PoleTele only $0.12$. There is no validation split and no early stopping; we report the final iterate. In none of the fits reported here did a length scale collapse or the noise reach its floor.

\subsection*{Gradient path}
The blocked and unblocked gradient implementations agree to $3\times10^{-5}$ in test RMSE. The unblocked path holds several $r \times n$ arrays in memory, while the blocked path processes the data in blocks of $20{,}000$ rows. We therefore use the blocked path by default and switch to the unblocked one only if a single timed evaluation shows it to be at least $20\%$ faster (\Cref{sec:experiments_new}). This margin does not guarantee that the unblocked path fits in memory; in the runs reported here it did.

\section{The predictive-mean solve}
\label{sec:solve}
\begin{figure}[h]
  \centering
  \includegraphics[width=\textwidth]{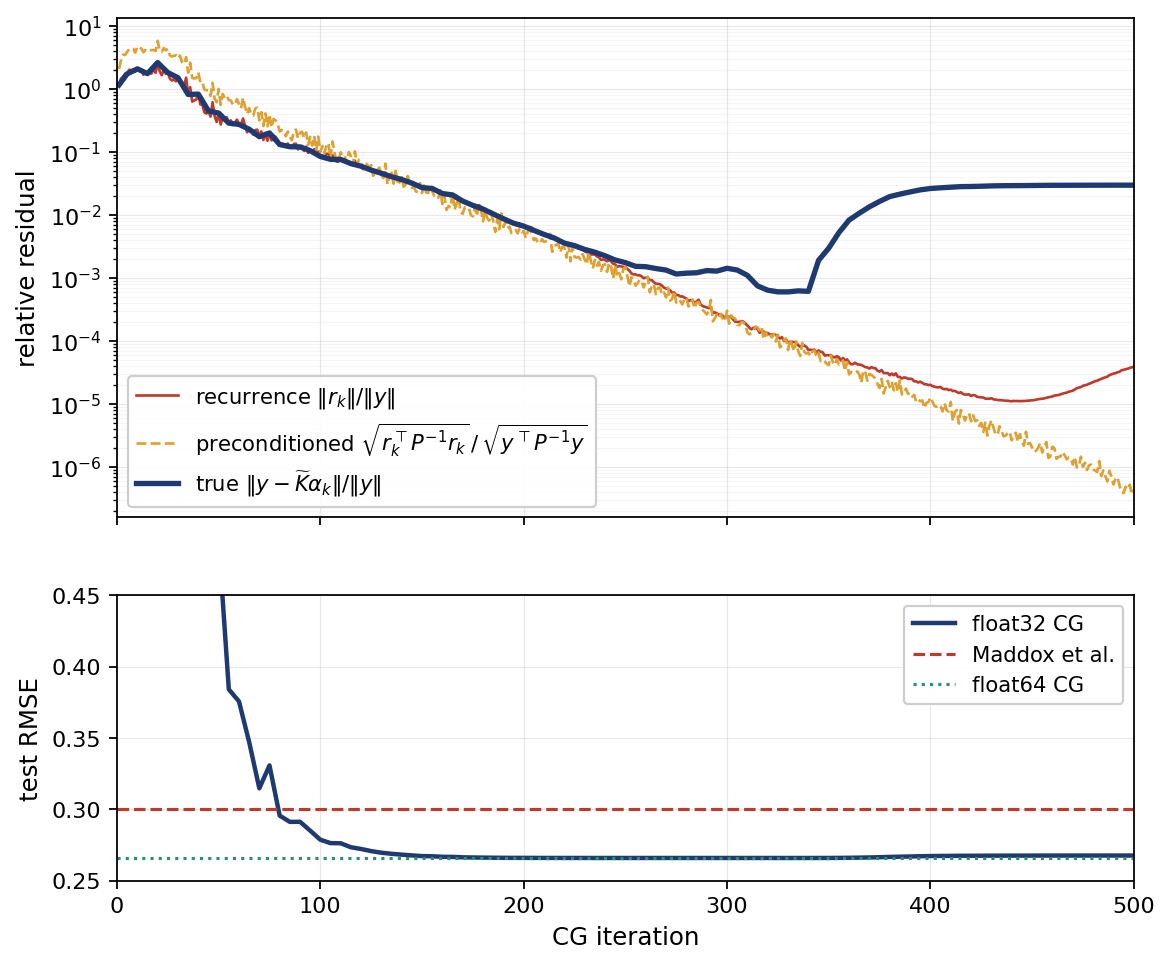}
  \caption{Residuals (top) and test RMSE (bottom) of fp32 CG with a rank-$200$ \AnchorNet\ preconditioner for Buzz with rank $250$.}
  \label{fig:buzz-resid}
\end{figure}

The exact predictive mean requires the solution of
\begin{equation}
  \widetilde{\bK}\,\boldsymbol{\alpha} = \by
  \label{eq:pm-solve}
\end{equation}
at the full training size. We use preconditioned CG with KeOps matrix--vector products, at most $200$ iterations, and a \AnchorNet\ preconditioner built from the first $\min(r, 200)$ landmarks of the model. \Cref{fig:buzz-resid} shows a fp32 solve on Buzz at rank $250$ over $500$ iterations. Up to about iteration $250$, the CG recurrence residual, the preconditioned residual norm, and the true residual $\lVert \by - \widetilde{\bK}\boldsymbol{\alpha}_k\rVert/\lVert \by\rVert$ agree. After that, the first two keep decreasing to $3.9\times10^{-5}$ and $4.4\times10^{-7}$, while the true residual stays at $3.0\times10^{-2}$. A stopping criterion based on either cheap quantity would therefore report convergence that has not happened. The preconditioned norm is particularly misleading, since $P^{-1}$ removes the residual component in the span of the landmarks, which matters most for prediction. The test RMSE is flat from iteration $200$ onward ($0.2659$ at iteration $200$, $0.2675$ at $500$), so more iterations do not help: the remaining error is a floor set by single precision.

\begin{figure}[h]
  \centering
  \includegraphics[width=\textwidth]{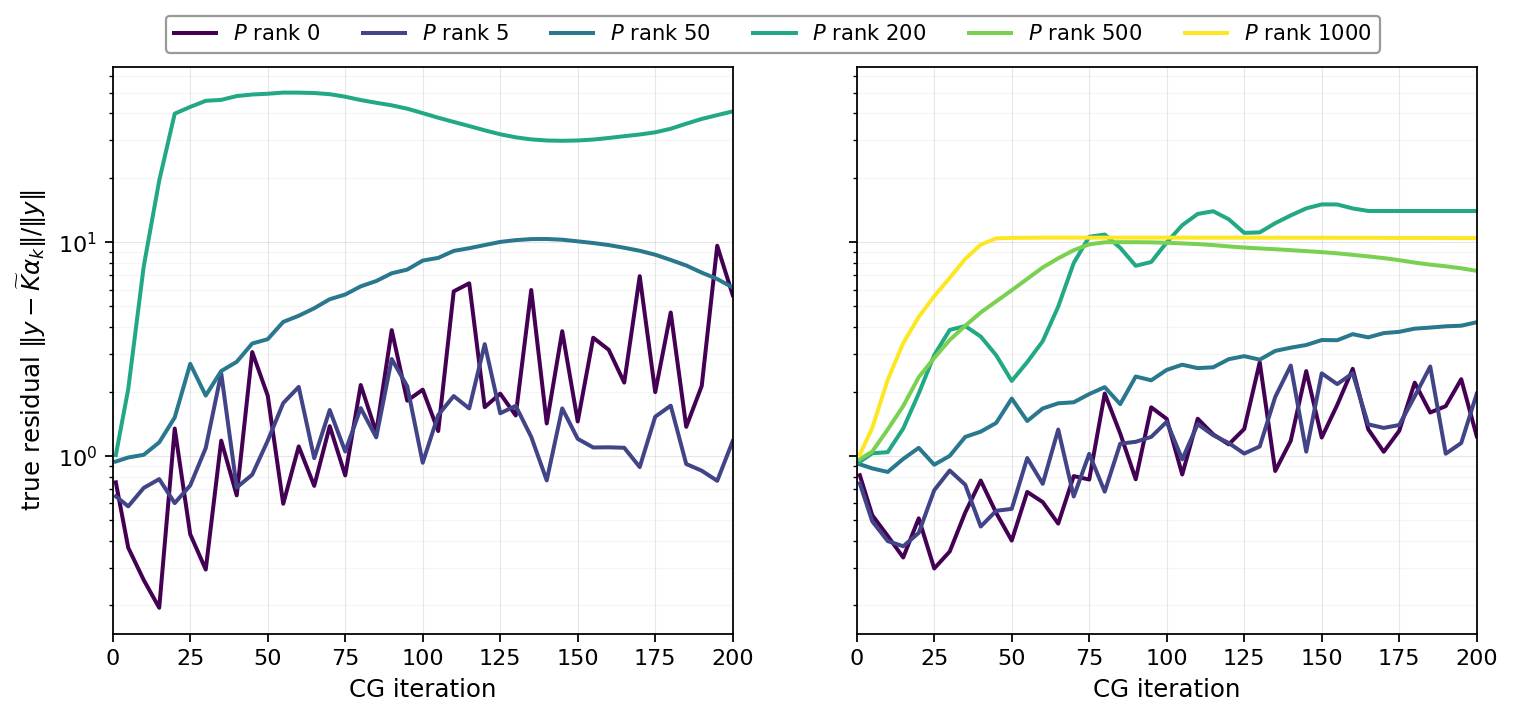}
  \caption{True residual of CG for different preconditioner ranks for HouseElectric in fp32, at model rank $500$ (left) and $3000$ (right). }
  \label{fig:he-precond}
\end{figure}

On HouseElectric, the fp32 solve does not converge for any preconditioner rank in $\{0, 5, 50, 200, 500, 1000\}$ below the model rank, at $r = 500$ and $r = 3000$ (\Cref{fig:he-precond}). The preconditioner here is built from the landmarks of the model, so its rank is at most the model rank. The true residual worsens as the preconditioner rank grows, and no preconditioner is the best setting. This is consistent with the factor $1/\sigma^2 \approx 400$ in $P^{-1}$, which amplifies its own rounding error.

\begin{figure}[h]
  \centering
  \includegraphics[width=\textwidth]{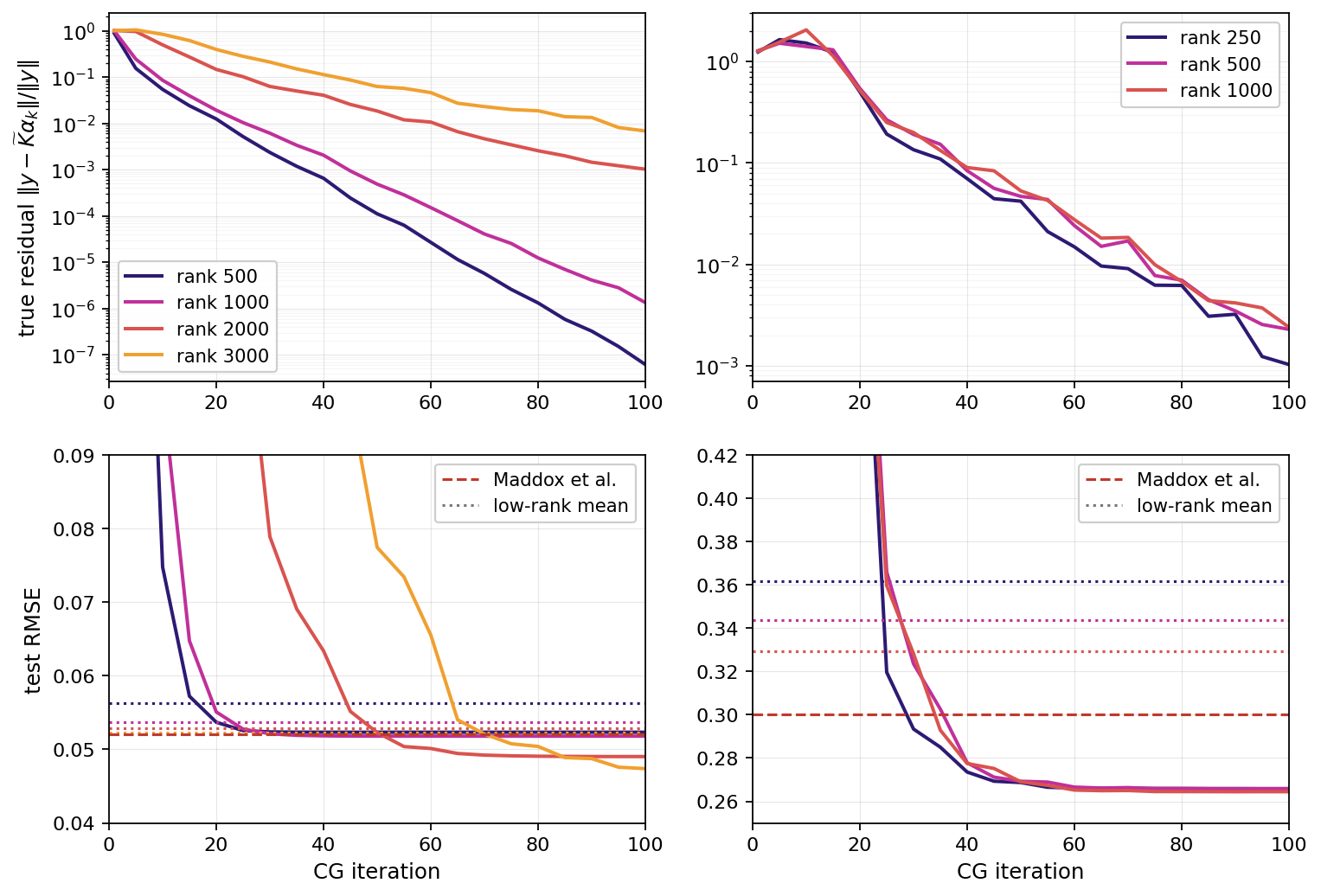}
  \caption{True residual (top) and test RMSE (bottom) of the fp64 solves on HouseElectric (left) and Buzz (right). Dotted lines show the low-rank predictive mean at the same rank, in the same color.}
  \label{fig:fp64}
\end{figure}

In fp64, with a rank-$1000$ preconditioner, the solve converges on both datasets (\Cref{fig:fp64}). On HouseElectric, the true residual reaches $6\times10^{-8}$ at rank $500$. These solves provide the reference values in Tables~\ref{tab:cost-3droad}--\ref{tab:cost-houseelectric}. On Buzz, the fp32 solve at rank $250$ gives $0.2723$, only $2.4\%$ above the fp64 value of $0.2659$, while the low-rank mean gives $0.3619$. There, the loss from single precision is much smaller than the loss from the low-rank predictive mean. On HouseElectric, the opposite holds: fp32 fails, and the low-rank mean is within $0.002$ to $0.006$ of the fp64 result from rank $200$ upward.

\end{document}

%% file: tables/summary.tex
\begin{table}[h]
  \centering
  \caption{Fit time comparison with \texttt{Maddox} \cite[Table 1]{maddox2022low}.}
  \label{tab:summary}
  \resizebox{\textwidth}{!}{%
  \begin{tabular}{lrrrrrrr}
    \toprule
    & & & \multicolumn{2}{c}{RMSE} & \multicolumn{3}{c}{time [s]} \\
    \cmidrule(lr){4-5}\cmidrule(lr){6-8}
    dataset & $(n, \, d)$ & $r$ & ANGP & \texttt{Maddox} & ANGP & \texttt{Maddox} & speed-up \\
    \midrule
    PoleTele & (13.5K, 26) & 100 & \textbf{0.1075\,$\pm$\,0.0085} & 0.117\,$\pm$\,0.004 & 2.2 & 58 & $26.6\times$ \\
    Elevators & (14.9K, 18) & 400 & \textbf{0.3607\,$\pm$\,0.0037} & 0.364\,$\pm$\,0.004 & 3.4 & 58 & $17.4\times$ \\
    Bike & (15.6K, 17) & 100 & \textbf{0.0707\,$\pm$\,0.0071} & 0.074\,$\pm$\,0.003 & 2.4 & 62 & $26.4\times$ \\
    Kin40K & (36K, 8) & 1600 & \textbf{0.0981\,$\pm$\,0.0018} & 0.099\,$\pm$\,0.001 & 17.4 & 66 & $3.8\times$ \\
    Protein & (41.2K, 9) & 1600 & \textbf{0.5406\,$\pm$\,0.0068} & 0.550\,$\pm$\,0.006 & 20.5 & 71 & $3.5\times$ \\
    3DRoad & (391.4K, 3) & 3000 & \textbf{0.1915\,$\pm$\,0.0040}$^{\ast}$ & 0.194\,$\pm$\,0.010 & 373 & 1,260 & $3.4\times$ \\
    Song & (463.8K, 90) & 3000 & \textbf{0.7605\,$\pm$\,0.0031} & 0.761\,$\pm$\,0.004 & 1,404 & 24,357 & $17.4\times$ \\
    Buzz & (524.9K, 77) & 50 & \textbf{0.2672\,$\pm$\,0.0056}$^{\ast}$ & 0.300\,$\pm$\,0.010 & 22.7 & 25,436 & $1,120.3\times$ \\
    HouseElectric & (1844.4K,11) & 800 & \textbf{0.0518\,$\pm$\,0.0003}$^{\ast}$ & 0.052\,$\pm$\,0.002 & 431 & 42,751 & $99.1\times$ \\
    \bottomrule
  \end{tabular}}
\end{table}

%% file: tables/rmse_vs_rank.tex
\begin{table}[h]
  \centering
  \small
  \caption{Test RMSE from the posterior mean solve against rank. The last row reports the speedup for hyperparameter fitting at the smallest rank that achieves at least the same accuracy as \texttt{Maddox}.
  % Fit time speed-up at the smallest such rank is reported in the last row.
  }
  \label{tab:rmse-rank}
  \resizebox{\textwidth}{!}{%
  \begin{tabular}{rrrrrrrrrr}
    \toprule
     & PoleTele & Elevators & Bike & Kin40K & Protein & 3DRoad$^{\ast}$ & Song & Buzz$^{\ast}$ & HouseElectric$^{\ast}$ \\
    $r$ & (0.117) & (0.364) & (0.074) & (0.099) & (0.550) & (0.194) & (0.761) & (0.300) & (0.052) \\
    \midrule
    50 & 0.1251 & 0.3806 & 0.0924 & 0.2663 & 0.7535 & 0.6319 & 0.7740 & \textbf{0.2672} & 0.0551 \\
    100 & \textbf{0.1075} & 0.3771 & \textbf{0.0707} & 0.2135 & 0.7086 & 0.5389 & 0.7710 & \textbf{0.2630} & 0.0549 \\
    200 & \textbf{0.1015} & 0.3667 & \textbf{0.0497} & 0.1764 & 0.6765 & 0.4944 & 0.7685 & \textbf{0.2655} & 0.0537 \\
    250 & \textbf{0.1000} & 0.3649 & \textbf{0.0480} & 0.1743 & 0.6539 & 0.4828 & 0.7677 & \textbf{0.2642} & 0.0534 \\
    400 & \textbf{0.1011} & \textbf{0.3607} & \textbf{0.0406} & 0.1371 & 0.6151 & 0.3991 & 0.7661 & \textbf{0.2648} & 0.0526 \\
    500 & \textbf{0.1016} & \textbf{0.3592} & \textbf{0.0389} & 0.1306 & 0.5952 & 0.3805 & 0.7652 & \textbf{0.2638} & 0.0522 \\
    800 & \textbf{0.1011} & \textbf{0.3565} & \textbf{0.0349} & 0.1106 & 0.5649 & 0.3015 & 0.7640 & \textbf{0.2622} & \textbf{0.0518} \\
    1000 & \textbf{0.1034} & \textbf{0.3555} & \textbf{0.0334} & 0.1051 & 0.5524 & 0.2871 & 0.7634 & \textbf{0.2628} & \textbf{0.0516} \\
    1600 & \textbf{0.1043} & \textbf{0.3545} & \textbf{0.0303} & \textbf{0.0981} & \textbf{0.5406} & 0.2502 & 0.7620 &  & \textbf{0.0506} \\
    2000 & \textbf{0.1055} & \textbf{0.3542} &  & \textbf{0.0931} & \textbf{0.5343} & 0.2316 & 0.7613 &  & \textbf{0.0490} \\
    3000 & \textbf{0.1076} & \textbf{0.3539} &  & \textbf{0.0835} & \textbf{0.5282} & \textbf{0.1915} & \textbf{0.7605} &  & \textbf{0.0467} \\
    4000 &  &  &  & \textbf{0.0759} & \textbf{0.5247} & \textbf{0.1753} & \textbf{0.7598} &  &  \\
    6000 &  &  &  & \textbf{0.0683} & \textbf{0.5216} & \textbf{0.1511} &  &  &  \\
    8000 &  &  &  & \textbf{0.0649} & \textbf{0.5207} &  &  &  &  \\
    \midrule
    \textit{speed-up} & $26.6\times$ & $17.4\times$ & $26.4\times$ & $3.8\times$ & $3.5\times$ & $3.4\times$ & $17.4\times$ & $1120.3\times$ & $99.1\times$ \\
    \bottomrule
  \end{tabular}}
\end{table}

%% file: tables/rmse_vs_rank_lr.tex
\begin{table}[t]
  \centering
  \small
  \caption{Test RMSE from the low-rank predictive mean for the fits of Table~\ref{tab:rmse-rank}.}
  \label{tab:rmse-rank-lr}
  \begin{adjustbox}{max width=\textwidth}
  \begin{tabular}{rrrrrrrrrr}
    \toprule
     & PoleTele & Elevators & Bike & Kin40K & Protein & 3DRoad & Song & Buzz & HouseElectric \\
    $r$ & (0.117) & (0.364) & (0.074) & (0.099) & (0.550) & (0.194) & (0.761) & (0.300) & (0.052) \\
    \midrule
    50 & 0.3595 & 0.3983 & 0.1812 & 0.6752 & 0.7912 & 0.6931 & 0.8333 & 0.4272 & 0.1012 \\
    100 & 0.3162 & 0.3903 & 0.1566 & 0.5743 & 0.7623 & 0.6320 & 0.8227 & 0.3796 & 0.0708 \\
    200 & 0.2744 & 0.3772 & 0.1114 & 0.4130 & 0.7325 & 0.5771 & 0.8145 & 0.3701 & 0.0596 \\
    250 & 0.2632 & 0.3742 & 0.1107 & 0.3757 & 0.7230 & 0.5641 & 0.8111 & 0.3619 & 0.0586 \\
    400 & 0.2438 & 0.3668 & 0.0932 & 0.3222 & 0.7004 & 0.5232 & 0.8069 & 0.3547 & 0.0566 \\
    500 & 0.2308 & 0.3646 & 0.0864 & 0.2965 & 0.6825 & 0.5026 & 0.8040 & 0.3442 & 0.0560 \\
    800 & 0.2010 & \textbf{0.3606} & \textbf{0.0646} & 0.2512 & 0.6565 & 0.4479 & 0.7988 & 0.3323 & 0.0542 \\
    1000 & 0.1882 & \textbf{0.3585} & \textbf{0.0687} & 0.2250 & 0.6437 & 0.4241 & 0.7966 & 0.3292 & 0.0536 \\
    1600 & 0.1698 & \textbf{0.3572} & \textbf{0.0517} & 0.1805 & 0.6180 & 0.3811 & 0.7922 &  & 0.0534 \\
    2000 & 0.1549 & \textbf{0.3561} &  & 0.1613 & 0.6094 & 0.3627 & 0.7904 &  & 0.0526 \\
    3000 & 0.1416 & \textbf{0.3549} &  & 0.1344 & 0.5892 & 0.3302 & 0.7862 &  & 0.0522 \\
    4000 &  &  &  & 0.1198 & 0.5805 & 0.3062 & 0.7838 &  &  \\
    6000 &  &  &  & 0.1013 & 0.5692 & 0.2716 &  &  &  \\
    8000 &  &  &  & \textbf{0.0915} & 0.5598 &  &  &  &  \\
    \bottomrule
  \end{tabular}
  \end{adjustbox}
\end{table}

%% file: tables/cost_3droad.tex
\begin{table}[t]
  \centering
  \caption{Test RMSE and prediction time of the three predictive means on 3DRoad. The reference RMSE is $0.194$.}
  \label{tab:cost-3droad}
  \begin{adjustbox}{max width=\textwidth}
  \begin{tabular}{rrrrrrrrrr}
    \toprule
    & & & & \multicolumn{2}{c}{low-rank mean} & \multicolumn{2}{c}{float32 CG} & \multicolumn{2}{c}{float64 CG} \\
    \cmidrule(lr){5-6}\cmidrule(lr){7-8}\cmidrule(lr){9-10}
    $r$ & eff.\ & fit [s] & $\times$ref & RMSE & time [s] & RMSE & time [s] & RMSE & time [s] \\
    \midrule
    50 & 49 & 4.10 & $307.3\times$ & 0.6931\,$\pm$\,0.0051 & 0.01 & 0.6319\,$\pm$\,0.0211 & 33 & 0.6319\,$\pm$\,0.0211 & 47.6 \\
    100 & 100 & 6.26 & $201.1\times$ & 0.6320\,$\pm$\,0.0069 & 0.01 & 0.5389\,$\pm$\,0.0110 & 33 & 0.5389\,$\pm$\,0.0110 & 94.1 \\
    200 & 198 & 10.98 & $114.7\times$ & 0.5771\,$\pm$\,0.0030 & 0.03 & 0.4944\,$\pm$\,0.0115 & 33 & 0.4944\,$\pm$\,0.0115 & 110 \\
    250 & 249 & 13.32 & $94.6\times$ & 0.5641\,$\pm$\,0.0018 & 0.03 & 0.4828\,$\pm$\,0.0112 & 33 & 0.4828\,$\pm$\,0.0113 & 128 \\
    400 & 397 & 22.54 & $55.9\times$ & 0.5232\,$\pm$\,0.0062 & 0.06 & 0.3992\,$\pm$\,0.0194 & 33 & 0.3991\,$\pm$\,0.0195 & 487 \\
    500 & 496 & 27.45 & $45.9\times$ & 0.5026\,$\pm$\,0.0033 & 0.07 & 0.3805\,$\pm$\,0.0102 & 33 & 0.3805\,$\pm$\,0.0101 & 553 \\
    800 & 795 & 50.96 & $24.7\times$ & 0.4479\,$\pm$\,0.0073 & 0.15 & 0.3015\,$\pm$\,0.0112 & 33 & 0.3015\,$\pm$\,0.0113 & 1,420 \\
    1000 & 989 & 67.23 & $18.7\times$ & 0.4241\,$\pm$\,0.0051 & 0.19 & 0.2876\,$\pm$\,0.0064 & 33 & 0.2871\,$\pm$\,0.0062 & 1,588 \\
    1600 & 1573 & 135.55 & $9.3\times$ & 0.3811\,$\pm$\,0.0056 & 0.39 & 0.2557\,$\pm$\,0.0077 & 33 & 0.2502\,$\pm$\,0.0065 & 2,091 \\
    2000 & 1959 & 192.69 & $6.5\times$ & 0.3627\,$\pm$\,0.0072 & 0.51 & 0.2371\,$\pm$\,0.0050 & 33 & 0.2316\,$\pm$\,0.0066 & 2,469 \\
    3000 & 2921 & 372.89 & $3.4\times$ & 0.3302\,$\pm$\,0.0038 & 0.90 & 0.2041\,$\pm$\,0.0052 & 33 & \textbf{0.1915\,$\pm$\,0.0040} & 2,950 \\
    4000 & 3869 & 743.44 & $1.7\times$ & 0.3062\,$\pm$\,0.0030 & 1.42 & 0.1953\,$\pm$\,0.0069 & 33 & \textbf{0.1753\,$\pm$\,0.0019} & 3,221 \\
    6000 & 5732 & 1,419.64 & $0.9\times$ & 0.2716\,$\pm$\,0.0040 & 2.91 & \textbf{0.1822\,$\pm$\,0.0057} & 33 & \textbf{0.1511\,$\pm$\,0.0020} & 3,556 \\
    \bottomrule
  \end{tabular}
  \end{adjustbox}
\end{table}

%% file: tables/cost_buzz.tex
\begin{table}[t]
  \centering
  \caption{Test RMSE and prediction time of the three predictive means on Buzz. The reference RMSE is $0.300$.}
  \label{tab:cost-buzz}
  \begin{adjustbox}{max width=\textwidth}
  \begin{tabular}{rrrrrrrrrr}
    \toprule
    & & & & \multicolumn{2}{c}{low-rank mean} & \multicolumn{2}{c}{float32 CG} & \multicolumn{2}{c}{float64 CG} \\
    \cmidrule(lr){5-6}\cmidrule(lr){7-8}\cmidrule(lr){9-10}
    $r$ & eff.\ & fit [s] & $\times$ref & RMSE & time [s] & RMSE & time [s] & RMSE & time [s] \\
    \midrule
    50 & 46 & 22.71 & $1120.3\times$ & 0.4272\,$\pm$\,0.0236 & 0.01 & \textbf{0.2783\,$\pm$\,0.0182} & 1,343 & \textbf{0.2672\,$\pm$\,0.0056} & 17,142 \\
    100 & 83 & 39.68 & $641.0\times$ & 0.3796\,$\pm$\,0.0168 & 0.02 & \textbf{0.2709\,$\pm$\,0.0158} & 1,341 & \textbf{0.2630\,$\pm$\,0.0029} & 19,022 \\
    200 & 151 & 74.94 & $339.4\times$ & 0.3701\,$\pm$\,0.0139 & 0.04 & 0.3021\,$\pm$\,0.0833 & 1,345 & \textbf{0.2655\,$\pm$\,0.0045} & 20,687 \\
    250 & 175 & 93.16 & $273.0\times$ & 0.3619\,$\pm$\,0.0100 & 0.05 & \textbf{0.2723\,$\pm$\,0.0057} & 1,337 & \textbf{0.2642\,$\pm$\,0.0047} & 20,396 \\
    400 & 244 & 149.77 & $169.8\times$ & 0.3547\,$\pm$\,0.0160 & 0.08 & 0.3047\,$\pm$\,0.0265 & 1,344 & \textbf{0.2648\,$\pm$\,0.0039} & 20,945 \\
    500 & 281 & 186.65 & $136.3\times$ & 0.3442\,$\pm$\,0.0071 & 0.10 & 0.3132\,$\pm$\,0.0282 & 1,341 & \textbf{0.2638\,$\pm$\,0.0040} & 20,949 \\
    800 & 393 & 308.35 & $82.5\times$ & 0.3323\,$\pm$\,0.0080 & 0.18 & 0.3563\,$\pm$\,0.0582 & 1,344 & \textbf{0.2622\,$\pm$\,0.0032} &  22,623 \\
    1000 & 461 & 389.44 & $65.3\times$ & 0.3292\,$\pm$\,0.0033 & 0.25 & 0.3414\,$\pm$\,0.0271 & 1,336 & \textbf{0.2628\,$\pm$\,0.0030} & 22,060 \\
    \bottomrule
  \end{tabular}
  \end{adjustbox}
\end{table}

%% file: tables/cost_houseelectric.tex
\begin{table}[t]
  \centering
  \caption{Test RMSE and prediction time of the three predictive means on HouseElectric. The reference RMSE is $0.052$.}
  \label{tab:cost-houseelectric}
  \begin{adjustbox}{max width=\textwidth}
  \begin{tabular}{rrrrrrrrrr}
    \toprule
    & & & & \multicolumn{2}{c}{low-rank mean} & \multicolumn{2}{c}{float32 CG} & \multicolumn{2}{c}{float64 CG} \\
    \cmidrule(lr){5-6}\cmidrule(lr){7-8}\cmidrule(lr){9-10}
    $r$ & eff.\ & fit [s] & $\times$ref & RMSE & time [s] & RMSE & time [s] & RMSE & time [s] \\
    \midrule
    50 & 50 & 24.17 & $1768.7\times$ & 0.1012\,$\pm$\,0.0225 & 0.03 & 1.0219\,$\pm$\,0.3712 & 1,522 & 0.0551\,$\pm$\,0.0006 & 12,353 \\
    100 & 90 & 44.45 & $961.8\times$ & 0.0708\,$\pm$\,0.0020 & 0.06 & 1.4842\,$\pm$\,1.1477 & 834 & 0.0549\,$\pm$\,0.0006 & 11,040 \\
    200 & 157 & 86.92 & $491.8\times$ & 0.0596\,$\pm$\,0.0002 & 0.14 & 7.0671\,$\pm$\,1.3902 & 223 & 0.0537\,$\pm$\,0.0001 & 13,651 \\
    250 & 181 & 108.54 & $393.9\times$ & 0.0586\,$\pm$\,0.0007 & 0.17 & 3.8081\,$\pm$\,2.2943 & 336 & 0.0534\,$\pm$\,0.0006 & 13,651 \\
    400 & 252 & 191.95 & $222.7\times$ & 0.0566\,$\pm$\,0.0004 & 0.33 & 4.1274\,$\pm$\,0.2744 & 272 & 0.0526\,$\pm$\,0.0007 & 17,570 \\
    500 & 286 & 238.73 & $179.1\times$ & 0.0560\,$\pm$\,0.0006 & 0.40 & 3.9374\,$\pm$\,0.5704 & 289 & 0.0522\,$\pm$\,0.0002 & 19,528 \\
    800 & 393 & 431.20 & $99.1\times$ & 0.0542\,$\pm$\,0.0004 & 0.74 & 4.6517\,$\pm$\,1.1249 & 305 & \textbf{0.0518\,$\pm$\,0.0003} & 24,752 \\
    1000 & 449 & 544.37 & $78.5\times$ & 0.0536\,$\pm$\,0.0003 & 0.91 & 4.2526\,$\pm$\,2.4126 & 272 & \textbf{0.0516\,$\pm$\,0.0002} & 26,128 \\
    1600 & 778 & 1,022.81 & $41.8\times$ & 0.0534\,$\pm$\,0.0003 & 1.73 & 1.3297\,$\pm$\,0.0791 & 535 & \textbf{0.0506\,$\pm$\,0.0003} & 42,254 \\
    2000 & 1011 & 1,368.83 & $31.2\times$ & 0.0526\,$\pm$\,0.0003 & 2.26 & 1.0387\,$\pm$\,0.2655 & 1,812 & \textbf{0.0490\,$\pm$\,0.0002} & 56,746 \\
    3000 & 1723 & 2,459.89 & $17.4\times$ & 0.0522\,$\pm$\,0.0003 & 4.13 & 0.9021\,$\pm$\,0.0515 & 1,811 & \textbf{0.0467\,$\pm$\,0.0003} & 86,133 \\
    \bottomrule
  \end{tabular}
  \end{adjustbox}
\end{table}